\documentclass[12pt]{article}

\usepackage{latexsym}
\usepackage{a4wide}
\usepackage{amscd}
\usepackage{graphics}
\usepackage{amsmath}
\usepackage{amssymb}
\usepackage{esint}
\usepackage{mathrsfs}
\usepackage{amsthm}
\usepackage{bbm}
\usepackage{color}
\usepackage{accents}
\usepackage{enumerate}
\usepackage[pdftex]{hyperref}
\usepackage{soul}

\usepackage{tikz-cd}

\newcommand{\N}{\mathbb{N}}                     
\newcommand{\Z}{\mathbb{Z}}                     
\newcommand{\R}{\mathbb{R}}                     
\newcommand{\C}{\mathbb{C}}                     
\newcommand{\p}{\partial}             		
\newcommand{\x}{\times}            		
\newcommand{\ind}{\mathrm{ind\,}}               
\newcommand{\CZ}{\mathrm{CZ}}               
\newcommand{\Crit}{\mathrm{Crit}\,}             
\newcommand{\SH}{\mathrm{SH}} 
\newcommand{\RFH}{\mathrm{RFH}}             
\newcommand{\FH}{\mathrm{FH}}             
\newcommand{\FC}{\mathrm{FC}}             
\renewcommand{\H}{\mathrm{H}}             
\newcommand{\ev}{\mathrm{ev}}               

\renewcommand{\j}{\mathfrak{j}}

\newcommand{\PP}{\mathcal{P}}    
     
\newcommand{\MM}{\mathcal{M}}     
\newcommand{\NN}{\mathcal{N}}

\newcommand{\A}{\mathbb{A}}

\renewcommand{\AA}{\mathcal{A}}            

\newtheorem{thm}{\sc Theorem}[section]               
\newtheorem*{thm*}{\sc Theorem}               
\newtheorem{cor}[thm]{\sc Corollary}        
\newtheorem*{cor*}{\sc Corollary}        
\newtheorem{lem}[thm]{\sc Lemma}            
\newtheorem{prop}[thm]{\sc Proposition}     
\newtheorem{rem}[thm]{\sc Remark}           
\newtheorem{conj}[thm]{\sc Conjecture}      

\DeclareFontFamily{U}{mathx}{\hyphenchar\font45}
\DeclareFontShape{U}{mathx}{m}{n}{
      <5> <6> <7> <8> <9> <10>
      <10.95> <12> <14.4> <17.28> <20.74> <24.88>
      mathx10
      }{}
\DeclareSymbolFont{mathx}{U}{mathx}{m}{n}
\DeclareFontSubstitution{U}{mathx}{m}{n}
\DeclareMathAccent{\widecheck}{0}{mathx}{"71}
\DeclareMathAccent{\wideparen}{0}{mathx}{"75}

\AtEndDocument{\bigskip{\footnotesize

  \addvspace{\medskipamount}
\noindent\textsc{Seoul National University, Department of Mathematical Sciences, Research Institute in Mathematics, 
	08826 Seoul, South Korea} \par  
 \noindent \textit{E-mail address}: \texttt{\href{mailto:jungsoo.kang@snu.ac.kr}{jungsoo.kang@snu.ac.kr}} \par
}}
\numberwithin{equation}{section}

\title{On the strong Arnold chord conjecture for convex contact forms} 

\author{Jungsoo Kang}

\date{}

\begin{document}
\maketitle

\begin{abstract}
The original Arnold chord conjecture states that every closed Legendrian submanifold of the standard contact sphere $S^{2n-1}$ admits a Reeb chord with distinct endpoints with respect to any contact form. In this paper, we prove  this conjecture for contact forms induced by strictly convex embeddings into $\R^{2n}$ under the assumption that minimal periodic Reeb orbits are of Morse-Bott type. We also provide a counterexample when the convexity condition is not satisfied.

\end{abstract}


\section{Introduction}

Arnold raised the following chord conjecture in his seminal paper \cite{Arn86}.
\begin{conj}\label{conj:chord}
	Every closed Legendrian submanifold $\Lambda$ of the standard contact sphere $(S^{2n-1},\xi_{\mathrm{std}})$ has a Reeb chord with respect to any contact form $\alpha$ supporting $\xi_\mathrm{std}$.
\end{conj}
An orbit $\gamma:[0,T]\to S^{2n-1}$ of the Reeb vector field associated with $\alpha$ satisfying $\gamma(0)\in\Lambda$ and $\gamma(T)\in\Lambda$ is called a Reeb chord of $\Lambda$ with action $T>0$. 
Conjecture \ref{conj:chord} was confirmed by Mohnke \cite{Moh01}. Partial results on this conjecture were made by Ginzburg-Givental, Abbas, and Cieliebak  \cite{Giv90,Abb99,Cie02}. Furthermore, Conjecture \ref{conj:chord} has been studied for more general classes of contact manifolds, see \cite{Abb04,HT11,HT13,Rit13,Zho20} as well as the literature mentioned above.

In fact, Arnold's original conjecture asks for the existence of a Reeb chord intersecting $\Lambda$ twice, or equivalently a Reeb chord $\gamma$ with distinct endpoints, i.e.~$\gamma(0)\neq\gamma(T)$. We call this the strong chord conjecture. As pointed out in \cite{Moh01}, the strong chord conjecture for arbitrary contact manifolds fails in general. Therefore, in this paper, we restrict our attention to the standard contact sphere case, the setting of Arnold's conjecture. There has not been much progress on the strong chord conjecture. It was proved by Ziltener \cite{Zil16} for the standard contact form on $(S^{2n-1},\xi_{\mathrm{std}})$, which corresponds to the unit sphere in $\R^{2n}$, see also \cite{CM18} for a quantitative result. The main goal of this paper is to settle the strong chord conjecture for contact forms induced by strictly convex embeddings of $S^{2n-1}$ into $\R^{2n}$ whose minimal periodic orbits are of Morse-Bott type. Moreover, we  provide a non-convex counterexample to the strong chord conjecture in Remark \ref{rem:intro}.(b). This counterexample was pointed out to the author by  Hutchings.

Let $S\subset \R^{2n}=\C^n$ be a starshaped domain. We mean by a domain a compact subset with smooth boundary and with nonempty interior. As $S$ is starshaped, the boundary $\p S$ of $S$ is transverse to all lines through the origin, and the restriction of the one-form $\lambda:=\sum_{i=1}^n\frac{1}{2}(x_idy_i-y_idx_i)$ to $\p S$ is a contact form. Every contact form on $(S^{2n-1},\xi_\mathrm{std})$ arises in this way. 
In this paper, we consider the boundary $\p C$ of a stictly convex domain $C\subset\R^{2n}$ with contact form $\alpha:=\lambda|_{\p C}$. By strict convexity, we mean that the Gauss-Kronecker curvature of $\partial C$ is positive definite everywhere. 
Let $\A_{\min}(\p C)$ denote the minimal action among the actions of all periodic Reeb orbits on $\p C$. We refer to periodic Reeb orbits with action $\A_{\min}(\partial C)$ as minimal periodic Reeb orbits. We also write $\A_{\min}(\Lambda)$ for the minimal action among the actions of all Reeb chords of a Legendrian submanifold $\Lambda$ in $\p C$.

\begin{thm}\label{thm:chord}
Let $C$ be a strictly convex domain in $\R^{2n}$. Assume that the set of minimal periodic Reeb orbits on $\partial C$ has a connected component of Morse-Bott type (see Section \ref{sec:MB} for the definition). 
	 Then every closed Legendrian submanifold $\Lambda$ of $\p C$ admits a Reeb chord with action strictly lower than $\A_{\min}(\p C)$, i.e.~$\A_{\min}(\Lambda)<\A_{\min}(\p C)$. In particular, this Reeb chord has distinct endpoints.
\end{thm}


Note that we do not assume fillability of Legendrian submanifolds, nor the Morse-Bott condition on their Reeb chords. Theorem \ref{thm:chord} applies, for example, to the unit sphere, as was proved in \cite{Zil16}. The Morse-Bott condition is weaker than the nondegeneracy condition, and thus the above theorem applies when there is a nondegenerate minimal periodic Reeb orbit on $\partial C$. 
Theorem \ref{thm:chord} proves the strong chord conjecture for the boundaries of strictly convex domains satisfying the Morse-Bott condition. Although we expect that the conjecture holds for the boundaries of arbitrary convex domains, we could not achieve this due to technical difficulties explained in Remark \ref{rem:arbitrary}. We also want to mention that our proof shows $\A_{\min}(\Lambda)<\A_{\min}(\partial C)$ for any convex domain $C$ and any closed Legendrian submanifold $\Lambda\subset\partial C$ that does not intersect any minimal periodic orbit, see Remark \ref{rem:Lambda}. The non-strict inequality $\A_{\min}(\Lambda)\leq \A_{\min}(\partial C)$ is true for any convex domain $C$  any closed Legendrian submanifold $\Lambda\subset\partial C$, see Proposition \ref{prop:SH} and below.

On the other hand, the strong chord conjecture and the non-strict inequality $\A_{\min}(\Lambda)\leq \A_{\min}(\partial S)$ fail for the boundaries of starshaped domains $S$ in general, see Remark \ref{rem:intro}. Therefore, in addition to various conjectures and results regarding convex domains in symplectic geometry, e.g.~\cite{HWZ98, Vit00,AMO08,AO08,ABHS15,AKP19,HK19,AK22,Iri22,GHV22,AB23}, one can interpret this as a new distinguished property of (strictly) convex domains in comparison to starshaped domains or general symplectic manifolds with contact boundaries. 

The proof of Theorem \ref{thm:chord} crucially  uses the strict convexity condition and does not seem to work under the dynamical convexity condition (see \eqref{eq:dyn_conv}), which is a weaker condition as shown in \cite{CE22}. 
The proof proceeds in three steps as outlined below. 

\begin{enumerate}[(1)]
	\item This step corresponds to Section \ref{sec:RFH_capacity}. We first consider the case that $C$ is strictly convex and  periodic Reeb orbits on $\p C$ are nondegenerate. In this case, applying results in \cite{AK22,Iri22}, we obtain a Floer cylinder with respect to a certain Hamiltonian function and a cylindrical almost complex structure in the symplectization $\R\times \p C$ such that its positive asymptotic orbit is a periodic Reeb orbit $x^-_{\min}$  with  action $\A_{\min}(\p C)$.
	\item The second step corresponds to  Section \ref{sec:pf_A}. We deform   a Floer cylinder in Step 1 to obtain a Floer curve $u$ defined on the closed disk with punctures, see Figure \ref{fig:limit}. It has exactly one positive interior puncture, at which $u$ converges to $x^-_{\min}$. It also has at least one negative boundary puncture, at which $u$ converges to a Reeb chord of $\Lambda$. Using strict convexity of $C$ again, we show that the $d\alpha$-energy of $u$ is always positive. By Stokes' theorem, the negative asymptotic Reeb chord of $u$ has action strictly lower than $x^-_{\min}$, and in turn the inequality $\A_{\min}(\Lambda)<\A_{\min}(\p C)$ follows.
	\item The last step corresponds to Section \ref{sec:MB}. Let $C$ be a strictly convex domain such that the set of minimal periodic Reeb orbits on $\partial C$ has a Morse-Bott component. We approximate $C$ by a sequence $(C_k)_{k\in\N}$ of strictly convex domains such that there is a unique nondegenerate minimal periodic orbit on each $\p C_k$. Applying the previous steps to $C_k$, we obtain a Floer curve $u_k$ (which was denoted by $u$ in Step 2) for each $k\in\N$. The sequence $(u_k)_{k\in\N}$ converges to a Floer curve $u$ which has exactly one interior positive puncture and at least one negative boundary puncture. The curve $u$ converges to a minimal periodic Reeb orbit on $\p C$ at the positive interior puncture and to a Reeb chord of $\Lambda$ at the negative boundary puncture. Due to the Morse-Bott condition, we know asymptotic behaviors of $u$ at the positive puncture, which allows us to obtain the positivity of the $d\alpha$-energy of $u$. As in Step 2, this gives rise to a Reeb chord with action strictly lower than $\A_{\min}(\partial C)$.    
\end{enumerate}

\begin{rem}\label{rem:intro}
The following examples show that Theorem \ref{thm:chord} is sharp.
\begin{enumerate}[(a)]
\item The following example shows that the inequality $\frac{\A_{\min}(\Lambda)}{\A_{\min}(\p C)}<1$ in Theorem \ref{thm:chord} is optimal. For a positive integer $a\in\N$, let 
\[
E_a=\Big\{(z_1,z_2)\in\C^2 \mid \frac{|z_1|^2}{a}+{|z_2|^2}\leq1\Big\},\quad L_a=\Big\{(z_1,z_2)\in\C^2\mid |z_1|^2=|z_2|^2=\frac{a}{a+1}\Big\}.
\]
The Reeb vector field is given by $\frac{2}{a}\frac{\p}{\p\theta_1}+2\frac{\p}{\p\theta_2}$ on  $L_a\subset\p E_a$, where $\theta_i=\arg(z_i)$. Thus $L_a$ is foliated by $a\pi$-periodic Reeb orbits. Every periodic Reeb orbit in $L_a$ intersects the Legendrian knot 
\[
\Lambda_a=\{(z_1,z_2)\in L_a \mid \theta_1+\theta_2=0\}
\] 
$a+1$ times, and we have
\[
\A_{\min}(\p E_a)=\pi,\qquad \A_{\min}(\Lambda_a)=\frac{a\pi}{a+1},\qquad \sup_{a\in\N} \frac{\A_{\min}(\Lambda_a)}{\A_{\min}(\p E_a)}=1.
\] 
\item There is a contact form $C^0$-close to the the standard one for which the strong chord conjecture does not hold. This counterexample was pointed out by Hutchings. We slightly perturb the line segment $\{(t,1-t)\in\R^2\mid t\in[0,1]\}$ near $(\frac{2}{3},\frac{1}{3})$ in such a way that the resulting curve $\eta$ still passes through $(\frac{2}{3},\frac{1}{3})$,  is tangent to the line $\{(t,t-\frac{1}{3})\in\R^2\mid t\in\R\}$ at $(\frac{2}{3},\frac{1}{3})$, and is transverse to all straight lines through the origin, see Figure \ref{fig:eta}. We consider  
\[
\p S_\eta=\big\{(z_1,z_2)\in\C^2\mid(|z_1|^2,|z_2|^2)\in\eta\big\},\quad L=\big\{(z_1,z_2)\in\C^2\mid(|z_1|^2,|z_2|^2)=(\tfrac{2}{3},\tfrac{1}{3})\big\}.
\]
Then $\p S_\eta$ is the boundary of a starshaped (toric) domain which is $C^0$-close to the unit ball in $\C^2$ but not convex. 
\begin{figure}[htb]
\centering
\includegraphics[width=0.35\textwidth,clip]{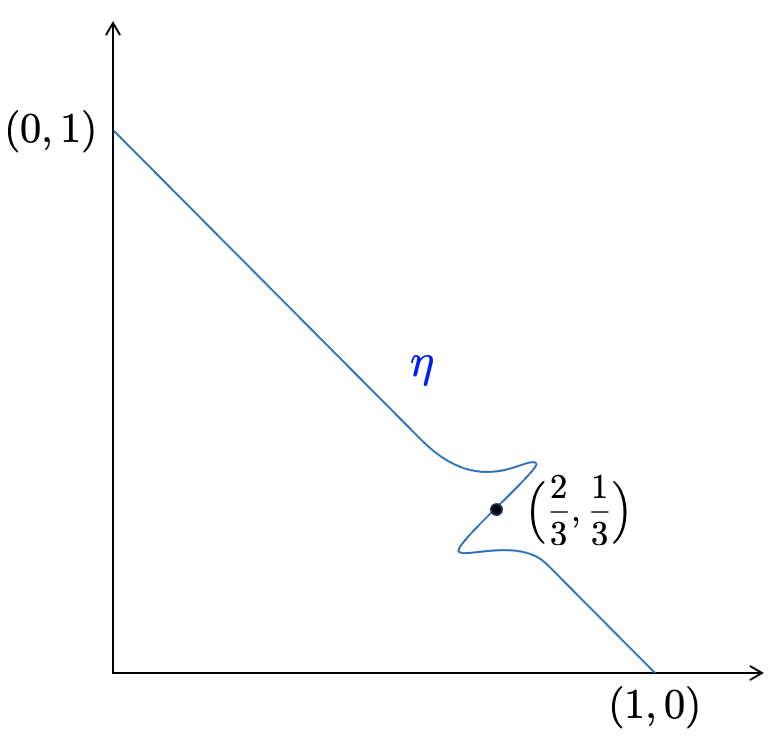}
\caption{A counterexample to the strong chord conjecture}\label{fig:eta}
\end{figure}
The Reeb vector field on $L$ is parallel to $\frac{\p}{\p\theta_1}-\frac{\p}{\p\theta_2}$, and every Reeb orbit in $L$ intersects the Legendrian knot 
\[
\Lambda=\{(z_1,z_2)\in L\mid 2\theta_1+\theta_2=0\}
\] 
exactly once. Thus $\Lambda$ does not admit a Reeb chord with distinct endpoints. We could even choose $\eta$ to pass through $(\frac{2}{3},\frac{1}{3}-\delta)$ for some small $\delta>0$ such that the tangent line at this point is parallel to the $x$-axis and has $(\frac{2}{3},\frac{1}{3}-\delta)$ as an isolated intersection point with $\eta$. Then the torus $\{(z_1,z_2)\in\C^2 \mid (|z_1|^2,|z_2|^2)=(\frac{2}{3},\frac{1}{3}-\delta)\}$ foliated by periodic orbits with minimal action $\A_{\min}(\partial S_\eta)=(\frac{1}{3}-\delta)\pi $ is of Morse-Bott type. Note that $\A_{\min}(\partial S_\eta)<\A_{\min}(\Lambda)=\pi$. One can also construct an analogous  counterexample in higher dimensions.

On the other hand, for a contact form sufficiently $C^1$-close to the standard one, every Legendrian torus has a Reeb chord with distinct endpoints, see \cite[Corollary 1.15]{CM18}. This tells that that  $C^0$-closeness in the above counterexample cannot be improved to $C^1$-closeness.

\item Let $B\subset\R^{2n}$ be the unit ball, and let $\Lambda$ be a closed Legendrian submanifold of $\p B$. By creating a little bump in a small Darboux open subset of $\p B$ away from $\Lambda$, we obtain a starshaped domain $S$ such that $\Lambda\subset \p S$ and $\p S$ contains a short periodic Reeb orbit. Thus,   $\frac{\A_{\min}(\Lambda)}{\A_{\min}(\p S)}$ can be made arbitrarily large.
\end{enumerate}
\end{rem}

It is considerably easier to establish the non-strict inequality $\A_{\min}(\Lambda)\leq \A_{\min}(\p C)$, as opposed to the strict inequality $\A_{\min}(\Lambda)< \A_{\min}(\p C)$, in Theorem \ref{thm:chord}. In fact, the non-strict inequality holds for more general class of contact manifolds as follows. We stress again that the strict inequality fails in general beyond convex domains as examples in \cite{Moh01} and Remark \ref{rem:intro} show.


\begin{prop}\label{prop:SH}
Let $Y$ be a closed smooth manifold with a contact form $\alpha$. If $(Y,\alpha)$ admits an exact filling $W$ with vanishing symplectic homology, then every closed Legendrian submanifold of $Y$ has a Reeb chord with action lower than or equal to $c_\SH(W)$.  
\end{prop}

The symplectic homology capacity $c_\SH$ is defined in Section \ref{sec:SH}. This proposition is not new and was proved by Zhou \cite{Zho20} using a different method. The action bound on a Reeb chord is not mentioned in \cite{Zho20} but should follow from his proof. 

Proposition \ref{prop:SH} implies  $\A_{\min}(\Lambda)\leq c_{\SH}(C)=\A_{\min}(\partial C)$ for every convex domain $C$ and every closed Legendrian submanifold $\Lambda\subset\partial C$, where the last equality follows from \cite{AK22,Iri22}.

We finally want to compare Proposition \ref{prop:SH} with Mohnke's theorem \cite{Moh01}. He proved that if a closed contact manifold $Y$ admits a subcritical Stein filling $W$, then every closed Legendrian submanifold of $Y$ admits a Reeb chord with action lower than or equal to the displacement energy of $W$, which is finite since a subcritical domain is displaceable in its completion, see  \cite{BK02}. Proposition \ref{prop:SH} generalizes his theorem  because a subcritical Stein domain has vanishing symplectic homology and the displacement energy is greater than or equal to $c_\SH$, see e.g.~\cite{Kan14}.

\paragraph{Acknowledgments} I am grateful to Alberto Abbondandolo for helpful comments and to Michael Hutchings for the counterexample in Remark \ref{rem:intro}.(b). I would like to thank the referees for their valuable comments and suggestions, which have helped improve the paper and correct an error. 
This work was supported by Samsung Science and Technology Foundation SSTF-BA1801-01 and National Research Foundation of Korea grant 2020R1A5A1016126 and RS-2023-00211186.

\section{Preliminaries}\label{sec:preliminary}
In this section, we collect some necessary properties of Rabinowitz Floer homology and symplectic homology. We refer to \cite{BO09,CF09,CFO10,CO18} for details of the construction of Floer homologies that we consider here.

A Liouville domain $(W,\lambda)$ is a compact smooth manifold $W$ of even dimension $2n$ with boundary $\partial W$ such that the exterior derivative $d\lambda$ of the one-form $\lambda$ on $W$ is symplectic and the Liouville vector field $Y$ characterized by $\iota_Yd\lambda=\lambda$ points outwards along $\partial W$. Using the flow $\phi_Y^t$ of $Y$, we can embed the negative half of the symplectization of $(\partial W,\alpha:=\lambda|_{\partial W})$ into $W$,
\[
(-\infty,0]\times \partial W\hookrightarrow W, \qquad (r,z)\mapsto \phi_Y^r(z),
\]
such that the pullback of $\lambda$ coincides with $e^r\alpha$. Attaching the other half of the symplectization by the flow $\phi_Y^t$, we obtain the completion $(\widehat W,\hat \lambda)$ of $(W,\lambda)$ defined by
\[
\widehat{W}:= W\cup_{\phi_Y}\big([0,\infty)\times\partial W\big),\qquad \hat\lambda= \begin{cases}
	\lambda & \text{on } W\\
	e^r\alpha & \text{on } [0,\infty)\times\partial W.
\end{cases}
\]
The contact form $\alpha=\lambda|_{\partial W}$ on $\partial W$ defines the Reeb vector field $R$ on $\partial W$ by $\iota_R d\alpha=0$ and $\iota_R\alpha=1$. We write $\phi_R^t$ for the flow of $R$. The action of a smooth $T$-periodic loop $\gamma:\R/T\Z\to\widehat  W$ is defined by 
\begin{equation}\label{eq:A}
\A(\gamma):=\int_0^T\gamma^*\hat\lambda.	
\end{equation}
If $\gamma$ is a Reeb orbit of $R$, then $\A(\gamma)$ coincides with its period $T$. We denote by $\mathrm{spec}(\partial W)$ the set of the actions of all periodic orbits of $R$. We set 
\[
\A_{\min}(\partial W):= \inf\mathrm{spec}(\partial W).
\]
Since $\partial W$ is compact, there exists a periodic Reeb orbit of $R$ whose action equals $\A_{\min}(\partial W)$. 
We assume that all periodic orbits $x$ of $R$ are nondegenerate, meaning that $d\phi_R^T(x(0))|_{\ker \alpha}$, the linearized return map restricted to the contact distribution $\ker \alpha$, does not have 1 as an eigenvalue, where $T>0$ is the period of $x$. We further assume that the first Chern class of the tangent bundle of $W$ vanishes on $\pi_2(W)$. This is necessary to have a well-defined index for periodic Reeb orbits contractible in $W$, but not indispensable for doing Floer theory.

\subsection{Rabinowitz Floer homology of $\partial W$}\label{sec:RFH}
\begin{figure}[htb]
\centering
\includegraphics[width=0.55\textwidth,clip]{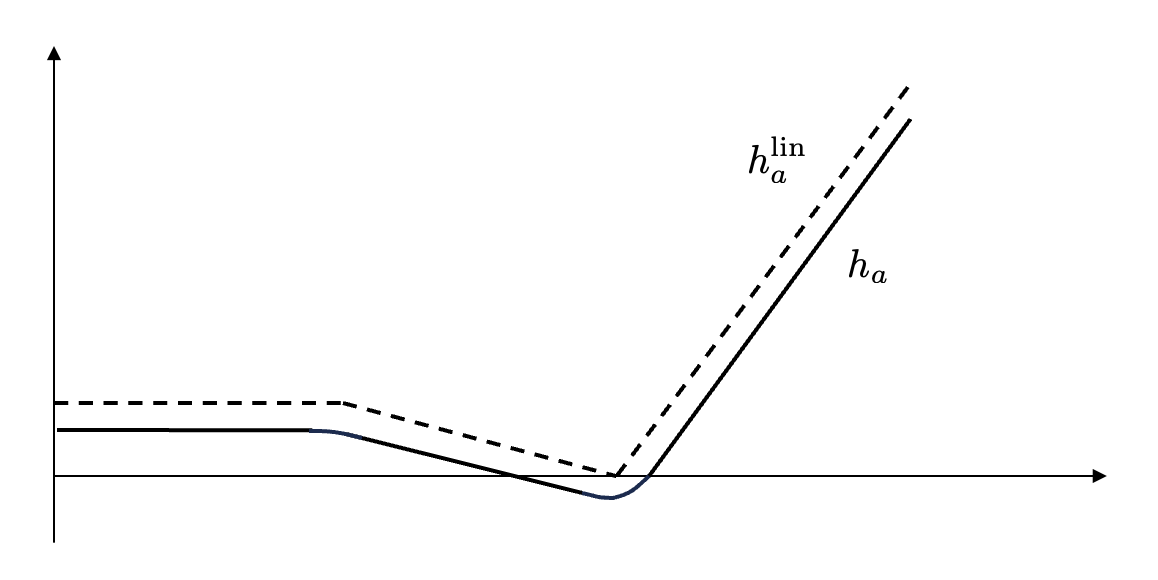}
\caption{graphs of $h_{a}^\mathrm{lin}$ and $h_a$}\label{fig:graph}
\end{figure}

To define the Rabinowitz Floer homology of $\partial W$ with action window $[0,\infty)$, we consider the following family of Hamiltonian functions. For $a\in\R\setminus \mathrm{spec}(\partial W)$, we consider a piecewise linear function $h_{a}^\mathrm{lin}:(0,+\infty)\to\R$ defined by 
\[
h_{a}^\mathrm{lin}(\tau) = 
\begin{cases}
	\epsilon/2 & \tau\in (0,1/2],\\[.5ex]
	-\epsilon\tau+\epsilon & \tau \in [1/2,1],\\[.5ex]
	a\tau-a & \tau \geq 1,
\end{cases}
\]
where $\epsilon=\epsilon_a\in(0,\A_{\min}(\p W))$ denotes a constant depending on $a$, which tends to zero as $a\to\infty$.
We take a smooth approximation $h_{a}$ of $h_{a}^\mathrm{lin}$, which is concave near $1/2$ and convex near $1$. We may further assume that $h_a$ attains its minimum at $1$ and $-\epsilon<\min h_a<0$, see Figure \ref{fig:graph}. 
We define $\overline H_{a}:\R\x \p W\to\R$ by $\overline  H_{a}(r,\cdot):= h_{a}(e^r)$. Remembering that $\R\x \p W$ is symplectically embedded inside  $\widehat{W}$, we define 
a smooth Hamiltonian $H_{a}:\widehat{W}\to\R$ by
\begin{equation}\label{eq:H_a}
	H_{a} =
\begin{cases}
	\overline{H}_{a} & \text{on } \R\times\partial W,\\[.5ex]
	\epsilon/2 & \text{on } \widehat{W}\setminus (\R\times\partial W).
\end{cases}
\end{equation}
For notational convenience, we write $h=h_a$ and $H=H_{a}$ when $a$ does not play any role. We define the Hamiltonian vector field $X_H$ by $-dH=\iota_{X_H}d\hat{\lambda}$. Note that $X_H=0$ on $ \widehat{W}\setminus (\R\times\partial W)$ and 
\begin{equation}\label{eq:X_H=R}
	X_H(r,\cdot)=h'(e^r)R \quad \textrm{on }\;\R\times\partial{W}.
\end{equation}
There are three types of 1-periodic orbits of $X_H$: 
\begin{enumerate}[(I)]
	\item constant 1-periodic orbits on $W\setminus([\log 1/2,0]\times\partial W)$,
	\item constant 1-periodic orbits on $r=0$,
	\item nonconstant 1-periodic orbits near $r=0$.
\end{enumerate}
If $x:S^1:=\R/\Z\to\R\times\partial W\subset \widehat W$ is an  orbit of type III, then by \eqref{eq:X_H=R} it is an  orientation-preversing reparametrization of a periodic Reeb orbit of $R$ with action 
\begin{equation}\label{eq:action_periodic}
\A(x)=h'(e^{r(x)}),	
\end{equation}
where $r(x)\in\R$ is the $\R$-coordinate of $x$. Note also that 1-periodic orbits of $X_H$ are exactly the critical points of the action functional 
\begin{equation}\label{eq:ham_action}
\AA_H:C^\infty(S^1,\widehat{W})\to\R,\qquad \AA_H(x):=\int_0^1x^*\hat{\lambda}-\int_0^1H(x)dt.
\end{equation}
For 1-periodic orbits $x_1$ and $x_2$ of type II or III, we have the following implication:
\begin{equation}\label{eq:implication}
\A(x_1)=\A(x_2)\quad \Longrightarrow \quad r(x_1)=r(x_2) \text{ and thus } \AA_H(x_1)=\AA_H(x_2).
\end{equation}
This follows from \eqref{eq:action_periodic} and the fact that $h$ is convex near $\tau=1$.  Furthermore, we have
\begin{equation}\label{eq:action_order}
\AA_H(\mathrm{I})<0< \AA_H(\mathrm{II})<\epsilon< \AA_H(\mathrm{III})	<a+\epsilon
\end{equation}
where the first inequality means that the Hamiltonian action of any type I orbit is lower than 0, and  likewise for other inequalities. 

Orbits of type II and III come in Morse-Bott families, and to remedy this issue we choose auxiliary Morse functions. More precisely, type II orbits form a Morse-Bott critical  manifold canonically diffeomorphic to $\partial W$. We choose a smooth Morse function $f_{\p W}:\p W\to \R$,
which has a unique maximum point $p_{\max}\in\p W$, and think of it as a function on the space of type II orbits. Moreover, since by our assumption all periodic orbits of $R$ with action lower than $a$ are nondegenerate, type III orbits form a Morse-Bott critical manifold diffeomorphic to the union of circles,
\begin{equation}\label{eq:S_x}
S_x:=\{\, x(\cdot+\theta) \mid \theta\in S^1 \,\}. 	
\end{equation}
We choose a perfect Morse function $f_x$ on $S_x$ and denote by $x^-$ and $x^+$ the minimum and maximum points of $f_x$, respectively.


\medskip

Now we construct a Floer chain complex using contractible 1-periodic orbits of $X_H$ of type II and III. From now on, we will only consider orbits contractible in $\widehat W$, omitting the adjective ``contractible''. For $e>\epsilon$, 
the $\Z_2$-vector space $\FC^{[0,e)}_*(H)$ is generated by all critical points $p$ of $f_{\p W}$ and $x^\pm$ of $f_x$ for every type III orbit $x$ with $\AA_H(x)<e$. It is graded as follows:
\[
\mu(p):=\ind_{f_{\p W}}(p)-n+1,\qquad \mu(x^-):=\CZ(x),\qquad \mu(x^+):=\CZ(x)+1
\]
where $\ind_{f_{\p W}}$ and $\CZ$ denote the Morse and the Conley-Zehnder indices, respectively. Recall that $\dim W=2n$. Note that if $e>a+\epsilon$, then all type III orbits of $X_H$ contribute to $\FC^{[0,e)}(H)$, see \eqref{eq:action_order}. In the case of $e=\epsilon$, all the generators of $\FC^{[0,\epsilon)}(H)$ are the critical points of $f_{\partial W}$. For any $e<e'$, we have a canonical inclusion 
\begin{equation}\label{eq:chain_incl}
	\FC^{[0,e)}_*(H) \hookrightarrow \FC^{[0,e')}_*(H).
\end{equation}

We choose a smooth $S^1$-family $J=(J_t)_{t\in S^1}$ of $d\hat\lambda$-compatible almost complex structures  on $\widehat{W}$ that is cylindrical at infinity. By being cylindrical at infinity we mean that there exists $r_+>0$ such that, on $(r_+,\infty)\times\partial W$, $J$ is independent of $t$, preserves $\ker\alpha$, is invariant under the translation $(r,\cdot)\to (r+c,\cdot)$, and satisfies $J\p_r=R$. We refer to smooth solutions $u:\R\x S^1\to \widehat{W}$ of 
\begin{equation}\label{eq:Floer}
\p_su+J_t(u)(\p_tu-X_H(u))=0	
\end{equation}
as Floer cylinders, where $(s,t)$ are coordinates on $\R\x S^1$. We denote by 
\begin{equation}\label{eq:evaluation}
\ev_-(u):=\lim_{s\to -\infty}u(s,\cdot), \qquad \ev_+(u):=\lim_{s\to +\infty}u(s,\cdot)
\end{equation}
provided the limits exist. In this case, the energy of $u$ is computed as
\begin{equation}\label{eq:energy0}
E(u):=\int_{\R\x S^1}|\p_su|^2dsdt=\AA_H(\ev_+(u))-\AA_H(\ev_-(u)),	
\end{equation}
where the norm $|\cdot|$ depending on $t\in S^1$ is induced by the metric $d\hat\lambda(\cdot,J_t\cdot)$. Note that the equality in \eqref{eq:energy0} still holds when $\ev_-(u)$ is replaced by $\displaystyle\lim_{k\to\infty}u(s_k^-,\cdot)$ provided this limit exists for some $s_k^-\to-\infty$, and similarly for $\ev_+(u)$. 
 This follows from Stokes' theorem and the fact that the Floer equation \eqref{eq:Floer} is the $L^2$-gradient flow equation of $\AA_H$. Moreover $E(u)=0$ yields $\ev_+(u)=\ev_-(u)=u(s,\cdot)$ for all $s\in\R$.

For $x_1,x_2\in\Crit\AA_H$ of type III, we denote by 
\[
\widehat{\MM}(S_{x_1},S_{x_2})=\widehat{\MM}(S_{x_1},S_{x_2};H,J)
\] the set of solutions $u$ of the Floer equation \eqref{eq:Floer} subject to the asymptotic condition $\ev_-(u)\in S_{x_2}$, $\ev_+(u)\in S_{x_1}$, and $\partial_su\to 0$ as $s\to\pm\infty$. We also denote 
\[
{\MM}(S_{x_1},S_{x_2}):=\widehat{\MM}(S_{x_1},S_{x_2})/\R,
\]
where the $\R$-action is given by $u(\cdot,\cdot)\mapsto u(s+\cdot,\cdot)$ for $s\in\R$. The evaluation maps $\ev_\pm$ descend to this quotient space, and we use the same notation.  
Let $\psi_{f_x}^t$ be the positive gradient flow of $f_x$ with respect to some metric on $S_x$. For $p\in\Crit f_x$, we write $W^u(p)$ and $W^s(p)$ for the unstable and stable manifolds of $\psi_{f_x}$ at $p$, respectively. For $\bar x,\underline x\in\Crit\AA_H$ of type III,  $p\in\Crit f_{\bar x}=\{\bar x^-,\bar x^+\}$ and $q\in\Crit f_{\underline x}$, we define 
\[
\MM(p,q):=\bigcup_{m\in\N\cup\{0\}}\MM_m(p,q)
\]
where $\MM_0(p,q):=(W^s(p)\cap W^u(q))/\R$ and, for $m\geq1$,
\[
\begin{split}
\MM_m(p,q):=\bigcup\;  & W^s(p)\x_{\ev_+}\MM(S_{\bar{x}},S_{x_1})\,_{\ev_-}\!\!\x_{\psi_{x_1}\circ\ev_+}\MM(S_{x_1},S_{x_2})\,_{\ev_-}\!\!\x \cdots \\[.5ex]
&\cdots \x_{\psi_{x_{m-1}}\circ\ev_+}\MM(S_{x_{m-1}},S_{\underline{x}})\,_{\ev_-}\!\!\x W^u(q)
\end{split}
\]
with the convention $x_0=\bar x$. The $\R$-action in $\MM_0$ is given by $z\mapsto \psi_{f_x}^s(z)$ for $s\in\R$, and the union in $\MM_m$ is over all possible $m-1$ type III orbits  $x_1,\dots,x_{m-1}$.

For $p\in\Crit f_x$ and $q\in\Crit f_{\partial W}$, the moduli space $\MM(p,q)$ is defined in the same way, but $W^u(q)$ should be understood as the unstable manifold of the positive gradient flow $\psi_{f_{\partial W}}$ of $f_{\partial W}$ at $q$. For $p,q\in\Crit f_{\partial W}$, we set $\MM(p,q)=(W^s(p)\cap W^u(q))/\R$, where the $\R$-action is given by $z\mapsto \psi_{f_{\partial W}}^s(z)$ for $s\in\R$. 

For a generic choice of $J$, $f_{x_i}$, and $f_{\partial{W}}$, the moduli space $\MM(p,q)$ is a smooth manifold of dimension $\mu(p)-\mu(q)-1$. If moreover $\mu(p)-\mu(q)=1$, $\MM(p,q)$ is compact, and we  denote by $n(p,q)\in\Z_2$ the parity of $\MM(p,q)$. Note that if $\MM(p,q)\neq \emptyset$, then $\AA_{H}(p)\geq \AA_H(q)$.

\medskip

We define the boundary operator for the Floer chain complex by 
\[
\p:\FC^{[0,e)}_*(H)\to \FC^{[0,e)}_{*-1}(H),\qquad p\mapsto\sum_{q} n(p,q)q,
\]
where the sum ranges over all generators $q$ with $\mu(p)-\mu(q)=1$.   
 We denote the homology of this chain complex by 
  \[
  \FH_*^{[0,e)}(H):=\H_*(\FC^{[0,e)}(H),\partial ).
  \]

Now, we reintroduce the subscript $a$ that was previously omitted from $H_a$. For any $a,b\in\R\setminus \mathrm{spec}(\partial W)$ with $a<b$, we take a homotopy $(H_s)_{s\in\R}$ from $H_b$ to $H_a$, which is  monotonically decreasing outside a compact subset,  and  consider the Floer equation in \eqref{eq:Floer} with $X_H$ replaced by $X_{H_s}$. Counting solutions of this equation together with gradient flow lines of $f_{\partial W}$ and $f_x$, we obtain  a continuation homomorphism 
\begin{equation}\label{eq:continuation}
c_{a,b}: \FH_*^{[0,e)}(H_a)\to \FH_*^{[0,e)}(H_b).	
\end{equation}
This family of vector spaces $\FH_*^{[0,e)}(H_a)$ together with homomorphisms $c_{a,b}$ forms a direct system. The limit of this system 
\[
\RFH_*^{[0,e)}(\partial W) := \lim_{a\to \infty}\FH_*^{[0,e)}(H_a)
\]
is referred to as the Rabinowitz Floer homology of $\partial W$ with action window $[0,e)$. In fact, if $\epsilon=\epsilon_e\in(0,\A_{\min}(\p W))$ is chosen sufficiently small, then, for any $e\leq a<b$, $c_{a,b}:\FH_*^{[0,e)}(H_a)\to \FH_*^{[0,e)}(H_b)$ is an isomorphism, see for instance \cite[Lemma 2.8]{Web06}. Therefore, we have 
\[
\RFH_*^{[0,e)}(\partial W)\cong \FH_*^{[0,e)}(H_a) \qquad \forall a\geq e. 
\] 
Moreover, the inclusion in \eqref{eq:chain_incl} induces homomorphisms
$\FH_*^{[0,e)}(H_a)\to \FH_*^{[0,e')}(H_a)$ for any $e<e'$,
and, by taking the limit over $a$, we obtain
\begin{equation}\label{eq:iota}
	\iota_{e,e'}:\RFH_*^{[0,e)}(\partial W)\to \RFH_*^{[0,e')}(\partial W).
\end{equation}
This family of homomorphisms also form a direct system, and we set
\[
\RFH_*^{[0,\infty)}(\partial W):= \lim_{e\to \infty}\RFH_*^{[0,e)}(\partial W).
\]
We remark that if one wants to define Rabinowitz Floer homology with full action window $(-\infty,\infty)$, then one should let $\epsilon$ in $h_a^{\mathrm{lin}}|_{[1/2,1]}=-\epsilon\tau+\epsilon$ go to infinity as $a\to\infty$.

\medskip

Next, we recall the Rabinowitz  Floer homology capacity from \cite[Definition A.1]{CO18}.  Note that $\FC^{[0,\epsilon)}(H)$ is isomorphic to the Morse chain complex of $f_{\p W}$ as  $\epsilon<\A_{\min}(\p C)$. 
  Thus, $\RFH^{[0,\epsilon)}_*(\partial W)$ is isomorphic to the singular homology $\H_*(\p W;\Z_2)$. Let us denote by $1_{\p W}$ the fundamental class of $\H_*(\p W;\Z_2)\cong \RFH^{[0,\epsilon)}_*(\partial W)$. Using the map in \eqref{eq:iota}, we define the Rabinowitz  Floer homology capacity of $\p W$ by
\[
c_\RFH(\partial W):=\inf\{e>0 \mid  \iota_{\epsilon,e}(1_{\p W})=0\}.
\]

\subsection{Symplectic homology of $W$}\label{sec:SH}

A goal of this section is to recall the definition of the symplectic homology capacity $c_\SH(W)$ and to prove that it coincides with the Rabinowitz Floer homology capacity $c_\RFH(\partial W)$ for Liouville domains $W$. To this end, we only need a specific diagram established in \cite{CFO10}. Therefore, we aim to keep the exposition of the construction of the symplectic homology $\SH_*(W)$ of $W$ as minimal as possible.

For $a\in\R\setminus \mathrm{spec}(\partial W)$,  we consider a piecewise linear function $k_{a}^\mathrm{lin}:(0,+\infty)\to\R$ defined by 
\[
k_{a}^\mathrm{lin}(\tau) = 
\begin{cases}
	\epsilon\tau-\epsilon & \tau \in (0,1],\\[.5ex]
	a\tau-a & \tau \geq 1,
\end{cases}
\]
where $\epsilon=\epsilon_a\in(0,\A_{\min}(\p W))$ tends to zero as $a\to\infty$.   
We take a smooth approximation $k_a$ of $k_{a}^\mathrm{lin}$ which is convex near $1$. 
Then we define a smooth Hamiltonian $K_a:\widehat{W}\to\R$ to have the following properties for a small $\delta>0$. 
\begin{enumerate}[(i)]
	\item On $\widehat{W}\setminus((-\delta,\infty)\times\partial W)$, $K_a$ is a $C^2$-small Morse function and takes values in $(-\epsilon,0)$.
	\item On $(-\delta,\infty)\times\partial W$, $K_a(r,\cdot)=k_a(e^r)$.
\end{enumerate}

For notational convenience, we omit the subscript $a$ from $k_a$ and $K_{a}$.  Note that a 1-periodic orbit of $X_{K}$ is either a constant orbit corresponding to a critical point of $K$ or a nonconstant orbit $x$ corresponding to a periodic Reeb orbit with action $\A(x)=k'(e^{r(x)})$, cf.~\eqref{eq:action_periodic}. 
Nonconstant orbits come in a Morse-Bott $S^1$-family. We again write $S_x$ for each component containing $x$, see \eqref{eq:S_x}, and choose a perfect Morse function $f_x:S_x\to\R$ with critical points $x^-$ and $x^+$.
 
For $e>\epsilon$, the $\Z_2$-vector space $\FC_*^{(-\infty,e)}(K)$ is generated by all constant orbits $p$ of $X_K$, which we identify with critical points of $K$, and $x^\pm\in\Crit f_x$ for nonconstant orbits $x$ of $X_K$ with $\AA_{K}(x)<e$. The action functional $\AA_K$ is the one introduced in \eqref{eq:ham_action} with $H$ replaced by $K$. It is graded as follows:
\[
\mu(p):=\ind_{-K}(p)-n,\qquad \mu(x^-):=\CZ(x),\qquad \mu(x^+):=\CZ(x)+1.
\]
The boundary operator $\partial: \FC_*^{(-\infty,e)}(K)\to\FC_{*-1}^{(-\infty,e)}(K)$ is defined in a similar way as in the previous section. We denote the homology of this chain complex by 
\[
\FH_*^{(-\infty,e)}(K):=\H_*(\FC_*^{(-\infty,e)}(K),\partial).
\]

Now we reintroduce the subscript $a$ that was omitted from $K_a$. As before for $a,b\in\R \setminus\mathrm{spec}(\partial W)$ with $a<b$, we may assume that $K_a<K_b$, and a monotone homotopy between them defines a continuation homomorphism 
\[
c_{a,b}:\FH_*^{(-\infty,e)}(K_a)\to\FH_*^{(-\infty,e)}(K_b).
\]
A family of such homomorphisms form a direct system, and the direct limit
\[
\SH^{(-\infty,e)}_*(W):= \lim_{a\to\infty}\FH_*^{(-\infty,e)}(K_a)
\]
is referred to as the symplectic homology of $W$ with action window $(-\infty,e)$. For any $e<e'$, the inclusion $\FC_*^{(-\infty,e)}(K)\hookrightarrow\FC_*^{(-\infty,e')}(K)$ induces a homomorphism
\[
\jmath_{e,e'}:\SH^{(-\infty,e)}_*(W)\to \SH^{(-\infty,e')}_*(W).
\]
The direct limit of this family of homomorphisms is called the symplectic homology of $W$ and denoted by 
\[
\SH_*(W):=\lim_{e\to\infty} \SH^{(-\infty,e)}_*(W).
\]

We note that $\FH_*^{(-\infty,\epsilon)}(K_a)$ and hence $\SH_*^{(-\infty,\epsilon)}(W)$ are  isomorphic to the Morse homology of $-K_a|_W$. The latter is isomorphic to the singular homology of the pair $(W,\partial W)$. We denote by $1_W$ the fundamental class of $\H_*(W,\p W;\Z_2)\cong \SH_*^{(-\infty,\epsilon)}(W)$. The symplectic homology capacity of $W$ is defined by
\[
c_\SH(W):=\inf\{e>0 \mid  \jmath_{\epsilon,e}(1_W)=0\}.
\]

Let $H_a$ be the Hamiltonian considered in Section \ref{sec:RFH}. We may assume that $K_a \leq H_a$. Then a monotone homotopy from $H_a$ to $K_a$ induces a homomorphism $\FH_*^{(-\infty,e)}(K_a)=\FH_*^{[0,e)}(K_a)\to \FH_*^{[0,e)}(H_a)$ like \eqref{eq:continuation}. Note that all the generators of $\FH_*^{(-\infty,e)}(K_a)$ have positive action values of $\AA_{K_a}$, and thus the action window $(-\infty,e)$ plays the same role as $[0,e)$. Taking the direct limit over $a$, we obtain 
\begin{equation}\label{eq:SH_RFH}
	\SH_*^{(-\infty,e)}(W)=\SH_*^{[0,e)}(W)\to \RFH_*^{[0,e)}(\partial W).
\end{equation}

\begin{lem}\label{lem:c_SH}
	It holds that $c_\SH(W)=c_\RFH(\partial W)$. 
\end{lem}
\begin{proof}

The homomorphism in \eqref{eq:SH_RFH} fits in to the following commutative diagram of long exact sequences established in \cite[Proposition 1.4]{CFO10}. 
\begin{equation}\label{eq:diagram}
\begin{tikzcd}[row sep=1.5em,column sep=1.3em]
\cdots \arrow{r} & \H_{2n}(W;\Z_2)\arrow{r} \arrow[equal]{d} & \SH^{(-\infty,a)}_n(W)  \arrow{r}{} & \RFH^{[0,a)}_n(\partial W) \arrow{r}{}  &  \cdots \,\,\, \\
\cdots \arrow{r} & \H_{2n}(W;\Z_2)\arrow{r} & \H_{2n}(W,\p W;\Z_2)  \arrow{u}{ \jmath_{\epsilon,e}} \arrow{r}{  }  & \H_{2n-1}(\partial W;\Z_2) \arrow{u}{\iota_{\epsilon,e}} \arrow{r}{}  & \cdots \,.
\end{tikzcd}
\end{equation}
The bottom line is induced by the inclusion $\partial W\hookrightarrow  W$, and $1_W$ is mapped to
$1_{\partial W}$. Since $\H_{2n}(W;\Z_2)=0$, the diagram yields $c_\SH(\partial W)= c_\RFH(W)$. 
\end{proof}

\subsection{RFH-capacity for convex domains}\label{sec:RFH_capacity}
From now on, we consider the situation that a Liouville domain $W$ is a strictly convex domain $C$ in $\R^{2n}$ with $\lambda=\sum_{i=1}^n\frac{1}{2}(x_idy_i-y_idx_i)$. We continue to write $\alpha:=\lambda|_{\partial C}$ and $R$ for the induced contact form and the associated Reeb vector field on $\partial C$, respectively. One important consequence of the strict convexity of $C$ is the dynamical convexity of $\partial C$, meaning that every periodic Reeb orbit $x$ on $\partial C$ has  
\begin{equation}\label{eq:dyn_conv}
\CZ(x)\geq n+1.	
\end{equation}
Moreover, every minimal periodic Reeb orbit $x_{\min}$, i.e.~$\A(x_{\min})=\A_{\min}(\partial C)$, has 
\begin{equation}\label{eq:CZ_min}
\CZ(x_{\min})= n+1.
\end{equation}
We refer to \cite{HWZ98,AK22} for the proofs.

The aim of this subsection is threefold. First, we observe  
the equality $c_\RFH(\partial C)=\A_{\min}(\partial C)$. Second, we explain that the Rabinowitz Floer homology $\RFH_*(\partial C)$ can be defined on the symplectization $\R\times\partial C$ without reference to the filling $C$. This second point will be needed later when we consider Legendrian submanifolds in $\partial C$ which do not necessarily admit Lagrangian fillings inside $C$. Finally, we show the existence of a special Floer cylinder, which will serve as the starting point of the proof of the main theorem.

\begin{cor}\label{cor:c_SH_convex}
	For every convex domain $C\subset\R^{2n}$, we have $c_\RFH(\partial C)=
	\A_{\min}(\partial C)$. 
\end{cor}
\begin{proof}
	It is proved in \cite{AK22,Iri22} that $c_\SH(C)=\A_{\min}(\partial C)$ holds for every  convex domain $C$. By Lemma \ref{lem:c_SH}, we conclude $ c_\RFH (\partial C)=c_\SH(C)=
	\A_{\min}(\partial C)$.
\end{proof}

Next, we briefly explain why $\RFH_*(\partial C)$ is well-defined on the symplectization $\R\times\partial C$. We use 
\begin{equation}\label{eq:H_a_symplectization}
H_a=\overline{H_a}: \R\times\partial C\to \R	
\end{equation}
given in \eqref{eq:H_a} and an $S^1$-family $J=(J_t)_{t\in S^1}$ of $d(e^r\alpha)$-compatible almost complex structures on $\R\times\partial C$, which are cylindrical outside $(r_-,r_+)\times\partial C$ for some $r_-<0<r_+$. The chain complex $(\FC^{[0,e)}_*(H_a),\partial)$ is defined in the same way as before, except now solutions $u:\R\times S^1\to \R\times \partial C$ of \eqref{eq:Floer} are involved in the definition of moduli spaces $\MM(p,q)$, rather than solutions in $\widehat{W}$. To ensure  $\partial\circ\partial=0$ in this case, we need to show that, for any moduli space $\MM(p,q)$ with $\mu(p)-\mu(q)\leq 2$, Floer cylinders in this space cannot escape to the negative end of the symplectization $\R\times\partial C$. This follows from the dynamical convexity of $\partial C$ mentioned in \eqref{eq:dyn_conv}, see \cite[Section 9.5]{CO18} for details. Furthermore, using a neck-stretching argument and the dynamical convexity of $\partial C$, one can see that the Floer homology defined in this way is isomorphic to the one we defined in the pervious section. We refer to the paragraph containing \eqref{eq:neck_stretch} for the neck-stretching argument in a comparable setting. 

We can make a similar argument for continuation homomorphisms in \eqref{eq:continuation}, defining them by counting solutions in $\R\times \partial C$. In this manner, we can construct $\RFH_*^{[0,e)}(\partial C)$ independently of the filling $C$, and it is isomorphic to the version defined in Section \ref{sec:RFH}.
\medskip

Now, we state a consequence of Corollary \ref{cor:c_SH_convex} that we will use in the proof of the main theorem. Before doing so, we observe that the unique maximum point $p_{\max}$ of the  Morse function $f_{\partial C}:\partial C\to\R$ (denoted $f_{\partial W}$ in Section \ref{sec:RFH}) represents
\[
[p_{\max}]=1_{\partial C}\in \H_{2n-1}(\p C;\Z_2)\cong \RFH^{[0,\epsilon)}_n(\partial C)\cong \FH_n^{[0,\epsilon)}(H_a)\quad \forall a>0.
\]
Note that $p_{\max}$ is the only generator of index $\CZ=n$ by \eqref{eq:dyn_conv}. Moreover, as mentioned in \eqref{eq:CZ_min}, every minimal periodic orbit $x_{\min}$ has $\CZ(x_{\min})=n+1$, and in turn 
\begin{equation}\label{eq:min_CZ}
\mu(x_{\min}^-)=n+1,\qquad \mu(x_{\min}^+)=n+2.
\end{equation}
\begin{lem}\label{lem:Floer_cylinder}
Let $H=H_a:\R\times \partial C\to\R$ for any $a>\A_{\min}(\partial C)$. Then there exists a type III orbit $x_{\min}$ of $X_H$ with $\A(x_{\min})=\A_{\min}(\partial C)$ such that $\MM(x_{\min}^-,p_{\max})$ has odd cardinality and consists of Floer cylinders $u$ with $\ev_+(u)=x_{\min}^-$ and $\ev_-(u)=p_{\max}$. 
  \end{lem}
  \begin{proof}
By the definition of $c_\RFH$ and Corollary \ref{cor:c_SH_convex}, for any $e>\A_{\min}(\partial C)$ sufficiently close to $\A_{\min}(\partial C)$, the map
 \[
 \iota_{\epsilon,e}:\RFH_n^{[0,\epsilon)}(\partial W)\cong \FH_n^{[0,\epsilon)}(H) \to \RFH_n^{[0,e)}(\partial W)\cong\FH_n^{[0,e)}(H)
 \]
 sends $1_{\partial W}$ to zero, i.e.~$p_{\max}$ is a boundary in $\FC_*^{[0,e)}(H)$. This together with \eqref{eq:min_CZ} yields that there is a minimal orbit $x_{\min}$ of type III such that $\MM(x_{\min}^-,p_{\max})$ has odd cardinality. Moreover, 
since $x_{\min}$ has minimal action, $\MM(x_{\min}^-,p_{\max})=\MM_1(x_{\min}^-,p_{\max})$ by \eqref{eq:energy0}, and every element $u$ in this moduli space satisfies 
\[
\qquad \ev_-(u)=p_{\max},\qquad \ev_+(u)\in W^s(x_{\min}^-)=\{x_{\min}^-\},
\]
where we use the fact that $x_{\min}^-$ is a minimum point of $f_{x_{\min}}$. 
This finishes the proof. 
  \end{proof}

\begin{rem}\label{rem:dyn_convex}
 Note that $x_{\min}^-$ in Lemma \ref{lem:Floer_cylinder} can be any generic reparametrization of $x_{\min}$ since we are free to choose a generic perfect Morse function $f_{x_{\min}}$ on $S_{x_{\min}}$. This fact will play a crucial role in the proof of the main theorem. 
	We stress that the strict convexity of $C$ is used twice in the proof of the lemma: $\ev_+(u) \in S_{x_{\min}}$ by Corollary \ref{cor:c_SH_convex} and $\ev_+(u) = x_{\min}^-$ by the dynamical convexity of $\p C$.
\end{rem}

\subsection{Floer curves from periodic orbits to chords}\label{sec:slit}

We will continue to focus on the case of strictly convex domains $C$, though the discussion below applies to general Liouville domains as well. In this section, we do not assume that periodic orbits of $R$ on $\p C$ are nondegenerate. Let $\Lambda$ be a closed Legendrian submanifold of $\p C$.
An orbit $y:[0,T]\to\partial C$ of the Reeb vector field $R$ with  $T>0$ is called a Reeb chord of $\Lambda$ if it satisfies $(y(0),y(T))\in\Lambda\times\Lambda$. The action of $y$ is defined by $\A(y):=\int_0^T y^*\alpha$, and it equals $T$. We denote by $\mathrm{spec}(\Lambda)$ the set of the actions of all Reeb chords of $\Lambda$. We set 
\[
\A_{\min}(\Lambda):= \inf\mathrm{spec}(\Lambda).
\] 
Like periodic orbits, Reeb chords are reparametrizations of  nonconstant Hamiltonian chords of $X_H$, where $H=H_a$ is the Hamiltonian considered in \eqref{eq:H_a_symplectization}. Here by a Hamiltonian chord we mean an orbit $y:([0,1],\{0,1\})\to (\R\times\partial C,\R\times\Lambda)$ of $X_H$.  Such Hamiltonian chords are exactly the critical points of the functional $\AA_H^\Lambda$ defined as follows.
Let $\PP_\Lambda$ be the space of smooth paths $y:[0,1]\to \R\x \p C$ with $y(i)\in\R\x\Lambda$ for $i=0,1$. We write $\AA_H^\Lambda$ for the action functional on $\PP_\Lambda$ having the same formula as $\AA_H$ in \eqref{eq:ham_action}, 
\begin{equation}\label{eq:A_lag}
\AA_H^\Lambda:\PP_\Lambda\to\R,\qquad \AA_H^\Lambda(y):=\int_0^1y^*(e^r\alpha)-\int_0^1H(y)dt.	
\end{equation}
If we take $a,\epsilon>0$ in the definition of $H=H_a$ to satisfy $a\notin\mathrm{spec}(\Lambda)$ and $\epsilon< \A_{\min}(\Lambda)$, then 
there are three types of critical points of $\AA_H^\Lambda$, as in the periodic case:
\begin{enumerate}[(I)]
	\item constant chords on $r<\log 1/2$,
	\item constant chords on $r=0$,
	\item nonconstant chords near $r=0$.
\end{enumerate}
The set of type II chords can be canonically identified with $\Lambda$. A type III chord $y$ corresponds to a Reeb chord of $\Lambda$ with action
\begin{equation}\label{eq:action_chord}
	\A(y)=\int_0^1 y^*\alpha=h'(e^{r(y)})
\end{equation}
via reparametrization, where $r(y)\in\R$ is the $\R$-coordinate of $y$. Type II chords form a Morse-Bott critical manifold for $\AA_H^\Lambda$, whereas type III chords need not be Morse-Bott.

Next, we consider curves that give rise to the 0-part of the closed-open map in Floer theory, see for example, \cite{Sei02, AS06a,Alb08,BC09,Abo10,Rit13} for a detailed account.   Following \cite{AS06a,AS10}, we work with the cylinder with a slit
\[
\Sigma := \R\times[0,1] \,/\, \R_{\geq0}\times \{0\} \sim \R_{\geq0}\times \{1\},
\]
see Figure \ref{fig:w} below. 
The surface $\Sigma$ carries a complex structure, which is regular everywhere, and a holomorphic coordinate $z=s+it$, which has exactly one singular point $(0,0)\sim (0,1)$. We refer to a smooth map $u:\Sigma\to \R\x\partial C$ solving the Floer equation \eqref{eq:Floer} on the interior of $\Sigma$ and satisfying the boundary condition
\[
u(s,t)\in\R\x \Lambda \qquad \forall (s,t)\in \p \Sigma = \R_{\leq 0}\x \{0,1\}\,/\,(0,0)\sim (0,1)
\]
as a Floer cylinder with a slit. In this case, $u$ is actually a solution of the Floer equation on the whole $\Sigma$, which makes sense even at $(0,0)\sim(0,1)$, see \cite{AS10} for details. 
As in \eqref{eq:evaluation}, we use the notation $\ev_\pm(u)=\displaystyle\lim_{s\to\pm\infty}u(s,\cdot)$ for the asymptotic orbits (if exist) of a Floer cylinder with a slit $u$. Note that $\ev_+(u)$ is a 1-periodic Hamiltonian orbit, and $\ev_-(u)$ is a Hamiltonian chord.

The energy $E(u)$ of a Floer  cylinder with a slit $u:(\Sigma,\partial\Sigma)\to(\R\x \p C,\R\times \Lambda)$ is defined by $E(u)=\int_\Sigma|\p_su|^2dsdt$, cf.~\eqref{eq:energy0}. If $E(u)$ is finite, then $u$ converges asymptotically to a 1-periodic orbit at the positive end and to a chord at the negative end. More precisely, there exist sequences $s^-_k\to-\infty$ and $s^+_k\to \infty$ such that 
\[
\ev_{(s_k^+)}(u):=\lim_{k\to\infty}u(s^+_k,\cdot)=x,\qquad \ev_{(s_k^-)}(u):=\lim_{k\to\infty}u(s^-_k,\cdot)=y
\]
for some 1-periodic orbit $x$ and chord $y$ of $X_H$. If $x$ is nondegenerate, then $\ev_+(u)=x$ holds, and likewise for $y$. By Stokes' theorem, the energy $E(u)$ of $u$ equals
\begin{equation}\label{eq:energy}
E(u)=\int_\Sigma|\p_su|^2dsdt=\AA_H(\ev_{(s_k^+)}(u))-\AA_H^\Lambda(\ev_{(s_k^-)}(u)).
\end{equation}

We will also use a Floer cylinder $u$ with a slit and also with negative punctures. By negative punctures, we mean the following. There is a finite subset $\Gamma:=\Gamma_{\mathrm{int}}\cup\Gamma_\p$ of $\Sigma$, where $\Gamma_{\mathrm{int}}\subset \Sigma\setminus\p\Sigma$ and $\Gamma_\p\subset\p \Sigma$. The map $u:\Sigma\setminus\Gamma\to\R\times\partial C$ is a solution of the Floer equation \eqref{eq:Floer} on the interior and has the following asymptotic behavior near $\Gamma$.  We pick a neighborhood $D_z$ of $z\in\Gamma_{\mathrm{int}}$ diffeomorphic to a closed disk and a biholomorphic map $\varphi_z:(-\infty,0]\x S^1\to D_z\setminus\{z\}$.  We write $\pi_\R$ and $\pi_{\p C}$ for the projections from $\R\x\p C$ to the first and the second factors, respectively. Then, there exists a sequence $\sigma^z_k\to-\infty$ such that
\begin{equation}\label{eq:asymptote}
	\lim_{\sigma\to-\infty}\pi_\R(u\circ \varphi_z(\sigma,\cdot))=-\infty,\qquad \ev_{(\sigma_k^z)}(u):=\lim_{k\to\infty}\pi_{\p C}(u\circ \varphi_z(\sigma_k^z,\cdot))=\gamma_z(T\,\cdot\,)
\end{equation}
	where $\gamma_z$ is a periodic Reeb orbit with action $T=\displaystyle\lim_{k\to \infty} \A(\pi_{\p C}(u\circ \varphi_z(\sigma_k^z,\cdot)))$. If $\gamma_z$ is nondegenerate, then $\displaystyle\lim_{\sigma\to-\infty}\pi_{\p C}(u\circ \varphi_z(\sigma,\cdot))=\gamma_z(T\,\cdot\,)$ in $C^\infty(S^1,\partial C)$ holds. A similar asymptotic property is required near $z\in\Gamma_\p$. We take a neighborhood $D_z$ of $z\in\Gamma_\p$ diffeomorphic to a closed half-disk and a biholomorphic map $\varphi_z:(-\infty,0]\x [0,1]\to D_z\setminus\{z\}$. Then, the properties in \eqref{eq:asymptote} hold, where now $\gamma_z$ is a Reeb chord of $\Lambda$. 
We call elements in $\Gamma_{\mathrm{int}}$ and $\Gamma_\p$ negative interior and boundary punctures, respectively.

If $u$ has finite $d\alpha$-energy, i.e.~$E_{d\alpha}(u)=\int_{\Sigma\setminus\Gamma}u^*d\alpha<\infty$, then $u$ has the aforementioned asymptotic properties, see \cite[Section 2.2]{Abb14}. There is also the notion of positive punctures, which we have not introduced, but in our case, all punctures will be negative. By Stokes' theorem, the $d\alpha$-energy of $u$ is computed as 
\begin{equation}\label{eq:dalpha_energy}
E_{d\alpha}(u)=\int_{\Sigma\setminus\Gamma}u^*d\alpha =  \A(\pi_{\p C}( \ev_{(s_k^+)}(u)))-\A(\pi_{\p C}( \ev_{(s_k^-)}(u)))-\sum_{z\in\Gamma}\A(\ev_{(\sigma^z_k)}(u)).
\end{equation}
We finally remark that if $J$ is cylindrical on the whole $\R\times\partial C$, then $X_H\in\ker d\alpha$, which yields
\begin{equation}\label{eq:positive_dalpha}
d\alpha(\p_su, \p_tu)= d\alpha(\p_tu-X_H(u), J(u)\p_tu)=d\alpha(\p_tu, J(u)\p_tu)\geq 0,
\end{equation}
and in turn $E_{d\alpha}(u)\geq 0$.


\section{Proofs of theorems}
\subsection{Proof of Theorem \ref{thm:chord} in the nondegenerate case}\label{sec:pf_A}
Let $C$ be a strictly convex domain in $\R^{2n}$, and let $\Lambda$ be a closed Legendrian submanifold in $\partial C$. We assume that all periodic Reeb orbits on $\p C$ with action equal to $\A_{\min}(\partial C)$ are nondegenerate. We take $H=H_a:\R\times S^1\to\R\times \partial C$ as in \eqref{eq:H_a_symplectization} to satisfy
\[
(0,a]\cap \mathrm{spec}(\p C)=\{\A_{\min}(\partial C)\},\qquad 0<\epsilon<\min\{\A_{\min}(\p C),\A_{\min}(\Lambda)\}.
\]
We take Morse functions $f_{\p C}$ and $f_{x}$ for every type III periodic orbit $x$ of $X_H$, which has $\A(x)=\A_{\min}(\partial C)$, as in the paragraph containing \eqref{eq:S_x}. Without loss of generality, we may assume that the maximum point $p_{\max}$ of $f_{\p C}$ and the minimum point $x^-$ of $f_{x}$ satisfy
\begin{equation}\label{eq:p_max}
p_{\max}\in\Lambda,\qquad p_{\max}\notin \R\x \mathrm{im\,}x,\qquad x^-(0)\notin\R\x\Lambda,	
\end{equation}
where $\mathrm{im\,}x$ denotes the image of $x$. 

\medskip

Let $\MM(x_{\min}^-,p_{\max})$ be the moduli space in Lemma \ref{lem:Floer_cylinder}. 
Let $c$ be the constant map from $\Sigma$ to $p_{\max}$, which is obviously a Floer cylinder with a slit. Since $p_{\max}$ is a Morse-Bott critical point of $\AA_H$ and also of $\AA_H^\Lambda$, the operator obtained by linearizing the Floer equation \eqref{eq:Floer} at $c$ is Fredholm of index zero. The index formula can be found in the references mentioned in Section \ref{sec:slit}. Moreover, as usual, this linearized operator at the trivial solution $c$ is automatically surjective, see \cite[Proposition 3.7]{AS06} for a proof in a comparable setting. We define 
\[
\MM:=\big\{([u],c)\mid [u]\in \MM(x_{\min}^-,p_{\max})\big\}.
\]
This space has odd cardinality by Lemma \ref{lem:Floer_cylinder}.

We consider the moduli space $\NN$, which consists of Floer cylinders with a slit $v:(\Sigma,\partial\Sigma)\to(\R\times\partial C,\R\times\Lambda)$ such that $\ev_+(v)=x^-_{\min}$,  $\ev_-(v)=p_{\max}$, and $v$ glued with any capping disk of $x^-_{\min}$ represents the trivial class in $\pi_2(\R\x C,\R\x\Lambda)$. We use a generic $S^1$-family $J_{\mathrm{reg}}$ of $d(e^r\alpha)$-compatible almost complex structures being cylindrical outside a compact set so that the moduli spaces $\MM$ and $\NN$ are cut out transversely. 
The index formula tells us that $\NN$ is one-dimensional. By the Floer-Gromov compactness and the gluing theorems, we know that $\NN$ is nonempty and $\MM$ is a boundary component of $\NN$.
\begin{figure}[htb]
\centering
\includegraphics[width=0.8\textwidth,clip]{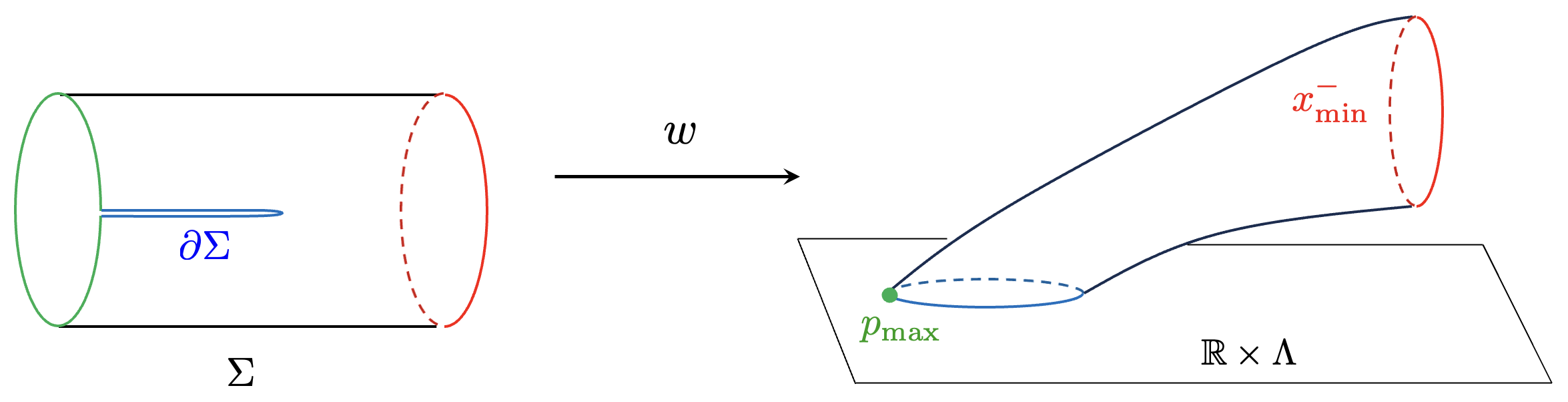}
\caption{Curves in $\NN$}\label{fig:w}
\end{figure}

\begin{rem}
	We remark that in an ideal situation, for example when $\Lambda$ admits an exact Lagrangian filling in $C$, the moduli space $\NN$ with its compactification in the sense of Floer-Gromov is used to obtain the commutative diagram in \eqref{eq:diagram2} below.
\end{rem}
   In our setting, we know that $\NN$ is compact up to breaking in the sense of Floer-Gromov and SFT, see e.g.~\cite{BO09a} and also \cite{BEHWZ03,CM05,Abb14}. Since the cardinality of $\MM$ is odd, there is a sequence $(v^m)_{m\in\N}$ in $\NN$ that does not converge, up to a subsequence, to an element of $\MM\cup\NN$. Note that, by the gluing theorem, it is not the case that a single element in $\MM$ compactifies a noncompact component of $\NN$ or serves as a boundary point of multiple components of $\NN$. We reparametrize $(v^m)$ in such a way that  it converges in the $C^\infty_{\rm loc}$-topology to a solution $u$ of the Floer equation \eqref{eq:Floer} such that $\ev_+(u)=x^-_{\min}$ and $u$ is not the  trivial cylinder at $x^-_{\min}$. We make the following observations on $u$. 
    \begin{enumerate}[(A)]
    	\item The domain of $u$ is either a punctured cylinder or a punctured cylinder with a slit. Interior and boundary punctures are due to bubbling of (punctured) pseudo-holomorphic spheres or disks, respectively. All punctures are negative by the maximum principle. The asymptotic behavior of $u$ near a puncture is described in Section \ref{sec:slit}.
    	\item  The domain of $u$ is a punctured cylinder if $v^m$ breaks along a 1-periodic orbit of $X_H$ in the sense of Floer theory. In this case, we claim  $\ev_-(u)=p_{\max}$ as a 1-periodic orbit. Indeed, $v^m$ cannot break along a type I periodic orbit since for any type I orbit $x_{\rm I}$ 
	\[
	\AA_H(x_{\rm I})<0<\AA_H(p_{\max})\leq \AA_H(v^m(s,\cdot))\leq \AA_H(x_{\min}^-) 
	\qquad \forall s\in \R.
	\]
	Note also that, for any type II and III orbits $x_{\rm II}$ and $x_{\rm III}$, we have $\AA_H(x_{\rm II})=\AA_H(p_{\max})$ and $\AA_H(x_{\rm III})=\AA_H(x_{\min}^-)$. This proves the claim. 
    	Moreover, since $(v^m)$ is chosen not to converge to an element in $\MM$,  $u$ must have a puncture.
    	\item  The domain of $u$ is a punctured cylinder with a slit if $v^m$ breaks along a chord of $X_H$ in the sense of Floer theory or does not break at all. 
    	\begin{itemize}
			\item[(C1)] If $v^m$ does not break, then $\ev_-(u)=p_{\max}$ as a chord, and there must be a puncture since $(v^m)$ is chosen not to converge to an element in $\NN$. 
			\item[(C2)] If $v^m$ breaks, then there exists a sequence $s^-_k\to-\infty$ such that $\displaystyle \ev_{(s^-_k)}(u)=\lim_{k\to\infty} u(s^-_k,\cdot)$ is a critical point of $\AA_H^\Lambda$. For a similar reason as in (B), $\ev_{(s^-_k)}$ must be a type III chord. 	
			\end{itemize}
    \end{enumerate}

   So far, we have worked with $J_{\mathrm{reg}}$, which may not be cylindrical on some compact subset due to a transversality issue. But, by a limiting argument, we may assume that the curve $u$ above solves the Floer equation \eqref{eq:Floer}  respect to  $J$ that is $d(e^r\alpha)$-compatible and cylindrical on the whole $\R\x \p C$. Indeed, we can take a sequence $(J^m_{\mathrm{reg}})_{m\in\N}$ converging to $J$ such that we have solutions $u^m$ (which was called $u$ above) of \eqref{eq:Floer} with respect to $J^m_{\mathrm{reg}}$ for each $m\in\N$. We reparametrize $(u^m)_{m\in\N}$ in such a way that, after passing to a subsequence, it converges to a map, which we call $u$ from now on, satisfying $\ev_+(u)=x^-_{\min}$ and  not being the trivial cylinder at $x^-_{\min}$. Then, this limit $u$ is a solution of   \eqref{eq:Floer} with respect to  $J$ and has the properties (A)-(C) mentioned above.

\begin{lem}\label{lem:domain_u}
	The domain of $u$ is a punctured cylinder with a slit. 
\end{lem}
\begin{proof}
In the proof, we use the contact action $\A$ defined in \eqref{eq:A}, rather than the Hamiltonian action $\AA_H$. 
Assume for contradiction that the domain of $u$ is $\R\x S^1\setminus\Gamma$, where $\Gamma$ is a finite subset of negative punctures. As in \eqref{eq:dalpha_energy} and \eqref{eq:positive_dalpha}, the $d\alpha$-energy of $u$ is computed as 
\[
E_{d\alpha}(u)=\int_{\R\x S^1\setminus\Gamma}u^*d\alpha =   \A(x_{\min}^-)-\A(\ev_{(s_k^-)}(u))-\sum_{z\in\Gamma}\A(\ev_{(\sigma^z_k)}(u))\geq 0,
\]
and $u^*d\alpha\geq 0$ holds pointwise. 
Thus, if $E_{d\alpha}(u)=0$, then $\p_tu$ and $\p_su$ are in the kernel of $d\alpha$, i.e.~$\p_tu,\p_su\in\R\langle R,{\p_r}\rangle$ at every point. Since $X_H$ and hence also $R$ are tangent to $\ev_+(u)=x^-_{\min}$, the image of $u$ belongs to the orbit cylinder of $x_{\min}^-$, i.e.~$\mathrm{im\,}u\subset\R\x\mathrm{im\,}x_{\min}^-$. 

By the observation (B) above, $\ev_-(u)=p_{\max}$ and $\Gamma\neq\emptyset$. Therefore we have $E_{d\alpha}(u)=0$ due to 
\[
\A(\ev_{(\sigma_k^z)}(u))\geq \A_{\min}(\p C)=\A(x_{\min}^-) \qquad \forall z\in\Gamma
\] 
and the above energy computation. This implies $\mathrm{im\,}u\subset\R\x\mathrm{im\,}x_{\min}^-$ as discussed above, but this contradicts  $\ev_-(u)=p_{\max}\notin\R\x \mathrm{im\,}x_{\min}^-$ in \eqref{eq:p_max}.
  This proves that the domain of $u$ cannot be $\R\x S^1\setminus\Gamma$. 
\end{proof}

\begin{figure}[htb]
\centering
\includegraphics[width=0.9\textwidth,clip]{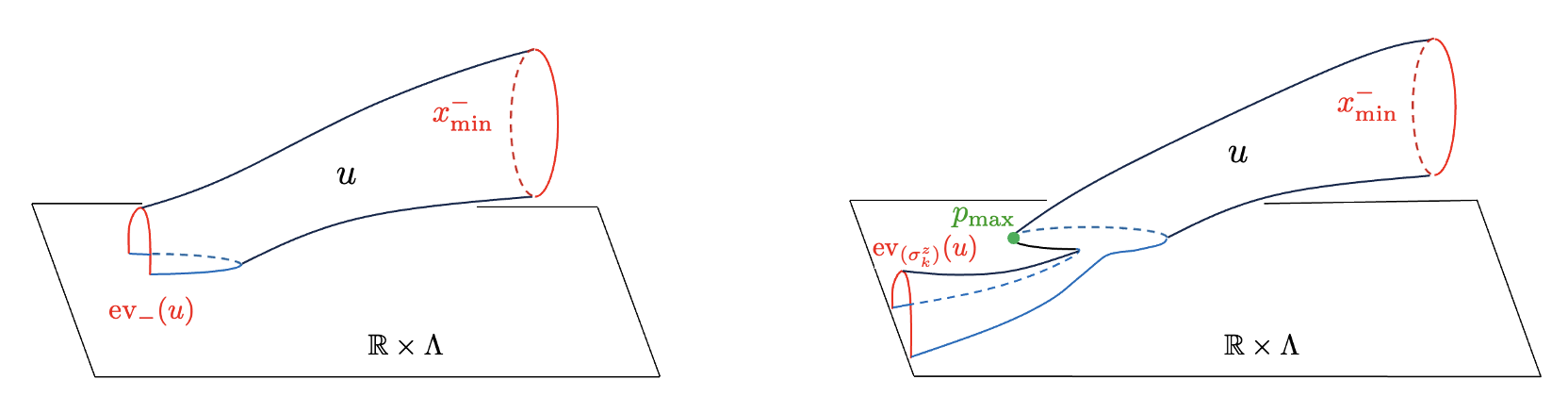}
\caption{Possible limit curves (i) and (ii)}\label{fig:limit}
\end{figure}

  Thus the map $u$ is defined on $\Sigma\setminus\Gamma$ with $\Gamma:=\Gamma_{\mathrm{int}}\cup\Gamma_\p$, where $\Gamma_{\mathrm{int}}$ and $\Gamma_\p$ are finite sets of negative punctures on $\Sigma\setminus\p\Sigma$ and $\p \Sigma$, respectively. By the observations (B) and (C) above, one of the following occurs, see also Figure \ref{fig:limit}.  
\begin{enumerate}[(i)]
	\item $\ev_{(s_k^-)}(u)$ is a type III chord of $X_H$.
	\item $\ev_-(u)=p_{\max}$ and $\Gamma=\Gamma_\mathrm{int}\cup\Gamma_{\p}$ is not empty.
	\end{enumerate}

\begin{lem}\label{lem:u_nontrivial}
	In case (i), we have  $\A(\ev_{(s_k^-)}(u))<\A(\ev_+(u))$. In case (ii), there exist $z\in\Gamma_\p$ and a sequence $\sigma_k^z\to-\infty$ such that $\A(\ev_{(\sigma_k^z)}(u))<\A(\ev_+(u))$.
\end{lem}
	\begin{proof}
	In case (i), if $\A(\ev_{(s_k^-)}(u))<\A(\ev_+(u))$ is not true, we have $\A(\ev_{(s_k^-)}(u))=\A(\ev_+(u))$ and $\Gamma=\emptyset$ due to $E_{d\alpha}(u)\geq0$. In particular, $E_{d\alpha}(u)=0$. Moreover, $E(u)=0$ and $\p_su=0$ by \eqref{eq:energy} since the implication in \eqref{eq:implication} also holds for a Floer cylinder with a slit because $\AA_H$ and $\AA^\Lambda_H$ have the same formula.  This leads to a contradiction that $\ev_{(s_k^-)}(u)(0)\in\R\x\Lambda$ while $\ev_+(u)(0)=x^-_{\min}(0)\notin\R\x\Lambda$ by \eqref{eq:p_max}. This proves $\A(\ev_{(s_k^-)}(u))<\A(\ev_+(u))$ as claimed. 

As observed in the proof of Lemma \ref{lem:domain_u}, $E_{d\alpha}(u)=0$ implies  $\mathrm{im\,}u\subset\R\x\mathrm{im\,}x_{\min}^-$.  Thus in case (ii) we have  $E_{d\alpha}(u)>0$  since $p_{\max}\notin\R\x \mathrm{im\,}x^-_{\min}$ by \eqref{eq:p_max}. Furthermore, $\Gamma_\mathrm{int}=\emptyset$ since otherwise $E_{d\alpha}(u)=0$ due to $\A(\ev_+(u))=\A_{\min}(\p C)$. Therefore there are  $z\in\Gamma_\p$ and a sequence $\sigma_k^z\to-\infty$ satisfying $\A(\ev_{(\sigma_k^z)}(u))<\A(\ev_+(u))$ by \eqref{eq:dalpha_energy}. 
	\end{proof}
	
 By Lemma \ref{lem:u_nontrivial}, $\ev_{(s_k^-)}(u)$ in case (i) or $\ev_{(\sigma_k^z)}$ for $z\in\Gamma_\p$ in case (ii) gives rise to a Reeb chord with action strictly lower than $\A_{\min}(\partial C)$. This proves Theorem \ref{thm:chord} under the assumption that all minimal periodic Reeb orbits are nondegenerate. \qed

	\begin{rem}
	One bad scenario that Lemma \ref{lem:u_nontrivial} excludes is that   $\ev_+(u)$ intersects $\R\times \Lambda$ and $\ev_{(s_k^-)}(u)(t)=\ev_+(u)(t)$ for all $t\in[0,1]$. That is, a periodic orbit $\ev_+(u)$ is also a chord of $\Lambda$ and coincides with $\ev_{(s_k^-)}(u)$. This is ruled out by the condition $\ev_+(u)(0)=x^-_{\min}(0)\notin\R\x\Lambda$ in \eqref{eq:p_max}. We stress again that the equality $\ev_+(u)(0)=x^-_{\min}(0)$ is a consequence of the dynamical convexity of $\partial C$, which is implied by the strict convexity of $C$, see Remark \ref{rem:dyn_convex}. Without the dynamical convexity of $\partial C$, it could be that $p_{\max}$ is the boundary of a generator of the form $x^+$, as opposed to $x^-$. In this case, we only know that there is  a Floer cylidner followed by a (possibly nontrivial) gradient flow line of $f_x$ connecting $p_{\max}$ and $x^+$, and in turn we cannot say $\ev_+(u)(0)\notin\R\times\Lambda$.  Thanks to the equality $\ev_+(u)(0)=x^-_{\min}(0)$ in our setting, even if $\ev_+(u)$ intersects $\R\times\Lambda$, we can conclude  $\A(\ev_{(s_k^-)}(u))<\A(\ev_+(u))$. 
	\end{rem}

\subsection{Proof of Theorem \ref{thm:chord} in the Morse-Bott case}\label{sec:MB}

Now, we assume that there is a connected component of the set of minimal periodic orbits of $R$ on $\partial C$ that is of Morse-Bott type. To be precise, there is a connected component 
\[
N\subset \{z\in\partial C\mid \phi_R^{\A_{\min}(\partial C)}(z)=z\},
\]
which is a closed submanifold of $\partial C$ such that $\ker d\alpha|_N$ has constant rank and $T_zN=\ker (d_z\phi_R^{\A_{\min}(\partial C)}-I)$. As in the previous section, let $H=H_a$ be the Hamiltonian associated to $\partial C$, and let $J$ be a $d(e^r\alpha)$-compatible  cylindrical almost complex structure on $\R\x\p C$. For $r_{\min}\in\R$ uniquely determined by $h'(e^{r_{\min}})=\A_{\min}(\partial C)$ (see \eqref{eq:action_periodic}), we set
\[
S_N:=\big\{x \mid x \text{ is a periodic orbit of $X_H$ with $x(0)\in \{r_{\min}\}\times N$}\big\}.
\]
It is a Morse-Bott critical manifold of the action functional $\AA_H$, and there is a natural diffeomorphism 
\[
\beta: S_N\to \{r_{\min}\}\times N\cong N,\qquad x\mapsto x(0).
\]

We perturb $\partial C$ to a strictly convex hypersurface with a nondegenerate minimal periodic orbit using a Morse function on the orbit space of $N$ as below, see \cite{Bou02} for details. To explain this, for a starshaped domain $S\subset\R^{2n}$, we consider a diffeomorphism $\rho:S^{2n-1}\to \p S$ where $S^{2n-1}$ is the unit sphere in $\R^{2n}$ and $\rho(z)$ is determined by $\rho(z)/|\rho(z)|=z$. Then $\alpha_S:=\rho^*(\lambda|_{\p S})=|\rho|^2(\lambda|_{S^{2n-1}})$ is a contact form on $(S^{2n-1},\xi_\mathrm{std})$. It is known that all contact forms on $(S^{2n-1},\xi_\mathrm{std})$ arise in this way.

	Let $\overline{N}:=N/(\R/\A_{\min}(\partial C)\Z)$ denote the  orbit space of $N$, and let $\pi:N\to \overline{N}$ be the associated circle bundle. We choose a negative smooth Morse function 
	\[
	\bar q:\overline{N}\to(-\infty,0),
	\] 
	which has a unique minimum point $z_{\min}\in \overline{N}$. 
	Consider the function $\bar q\circ \pi\circ \rho$ defined on $\rho^{-1}(N)\subset S^{2n-1}$. We extend it to $S^{2n-1}$ using a cutoff function in a tubular neighborhood of $\rho^{-1}(N)$.
	We denote the resulting function by $q:S^{2n-1}\to \R$. This function gives rise to a sequence of contact forms $\alpha_{C,k}:=(1+\frac{1}{k} q)\alpha_C$ for $k\in\N$, where $\alpha_C$ is the contact form on $S^{2n-1}$ induced by $C$. We denote by $C_k$ the starshaped domain   corresponding to $\alpha_{C,k}$ via the map $\rho$ above, i.e.~$\alpha_k:=\lambda|_{\partial C_{k}}=(\rho^{-1})^*\alpha_{C,k}$ holds. For a sufficiently large $k\in\N$,  the following properties hold.
	
	\begin{enumerate}[-]
		\item The domain $C_k$ is strictly convex. 
		\item There is a unique periodic Reeb orbit  $\gamma_{k,\min}$ on $\partial C_k$ with minimal action $\A(\gamma_{k,\min})=\A_{\min}(\partial C_{k})$, and it is nondegenerate. As $k\to\infty$, $\gamma_{k,\min}$ converges to a periodic Reeb orbit $\gamma_{\min}$ on $\partial C$ satisfying $\mathrm{im\,}\gamma_{\min}=\pi^{-1}(z_{\min})$, i.e.~$\gamma_{k,\min}(\A_{\min}(\partial C_{k})t)\to \gamma_{\min}(\A_{\min}(\partial C) t)$ in $C^\infty(S^1,\R^{2n})$. 
	\end{enumerate} 
	
	Let $x_{\min}$ denote a type III periodic orbit of $H$, corresponding to $\gamma_{\min}$ via reparametrization. Note that $x_{\min}$ is a minimum point of $\bar q\circ \pi\circ \beta: S_N\to(-\infty,0)$. We may assume that $x_{\min}(0)\notin\R\times\Lambda$.

	Let $H_k:\R\x\p C_k\to\R$ be a Hamiltonian associated to $\partial C_k$ defined in the same way as $H$. We may assume that $H_k$ converges to $H$ in the $C^\infty_{\rm loc}$-topology as $k\to\infty$. Let $x_{k,\min}$ be a type III periodic orbit of $H_k$, corresponding to $\gamma_{k,\min}$. We choose a perfect Morse function on $S_{x_{k,\min}}$ such that its minimum point $x_{k,\min}^-$ converges to $x_{\min}$ in $C^\infty(S^1,\R^{2n})$. Then $x_{k,\min}^-(0)\notin\R\times\Lambda_k$ for sufficiently large $k$, where $\Lambda_k:=\p C_k\cap(\R\x\Lambda)$. As usual, let $f_{\partial C_k}$ be a Morse function on $\partial C_k$ with a unique maximum point $p_{k,\max}$. As in \eqref{eq:p_max}, we may assume that $p_{k,\max}\in\Lambda_k\setminus( \R\times \mathrm{im\,} x_{k,\min})$, and additionally that $p_{k,\max}$ converges to $p_{\max}\in\Lambda\setminus( \R\times \mathrm{im\,} x_{\min})$ as $k\to\infty$.

We also choose a sequence $(J_k)_{k\in\N}$ of cylindrical $d(e^r\alpha_k)$-compatible almost complex structures on $\R\x\p C_k$ converging to $J$. 
	Note that all $(\R\x\p C_k,d(e^r\alpha_k))$ are symplectomorphic to $(\R^{2n}\setminus\{0\},d\lambda)$ but the forms of $H_k$ and $J_k$ depend on the splitting $\R\x\p C_k$. The proof in Section \ref{sec:pf_A} yields a sequence $(u_k)_{k\in\N}$ of Floer cylinders with a slit and possibly with negative boundary punctures in $(\R\x\p C_k,\R\x \Lambda_k)$ with respect to $(H_k,J_k)$ such that $\ev_+(u_k)=x^-_{k,\min}$ and either case (i) or case (ii) (see above Lemma \ref{lem:u_nontrivial}) is true for $u_k$.  After passing to a subsequence, a suitable reparametrization of $u_k$ converges to  a Floer cylinder with a slit and possibly with punctures $u$ that is of either type (i) or (ii) and satisfies the asymptotic condition 
\begin{equation}\label{eq:asymptotic_condition}
\ev_+(u)\in W^s(x_{\min}),
\end{equation}
	where $W^s(x_{\min})$ is the stable manifold of the gradient flow of $\bar{q}\circ\pi\circ\beta:S_N\to(-\infty,0)$. 
	This is the top layer of the limit of $u_k$ in the sense of SFT in the Morse-Bott case, see \cite[Proposition 4.16]{Bou02} and \cite[Section 4.2]{BO09}.  Note that we cannot apply                                                                                                                                                                                                                                                                                                                                                                                                      the SFT compactness theorem to lower layers due to the lack of the nondegeneracy (or the Morse-Bott property) of Reeb chords of $\Lambda$ but this is not necessary for our purpose. Since $x_{\min}$ is a minimum point of $\bar{q}\circ\pi\circ\beta$, the condition \eqref{eq:asymptotic_condition} reads $\ev_+(u)=x_{\min}$. Since $x_{\min}(0)\notin\R\times \Lambda$ and $p_{\max}\notin \R\times \mathrm{im\,} x_{\min}$, Lemma \ref{lem:u_nontrivial} applies to $u$ and we conclude that there is a Reeb chord with action strictly lower than $\A_{\min}(\partial C)$. This completes the proof of Theorem \ref{thm:chord}.  \qed

	\begin{rem}\label{rem:Lambda}
		Let $C\subset \R^{2n}$ be a convex domain, not necessarily strictly convex. We do not assume that the set of minimal periodic Reeb orbits on $\partial C$ has a Morse-Bott component. Instead, let $\Lambda$ be a closed Legendrian submanifold of $\partial C$, which does not intersect any minimal periodic Reeb orbit. Then $\A_{\min}(\Lambda)<\A_{\min}(\partial C)$ holds, and thus $\Lambda$ admits a Reeb chord with distinct end points. We outline the proof of this assertion. 
		
		We can approximate $C$ by a sequence $(C_k)_{k\in\N}$ of strictly convex domains such that there is a unique minimal periodic Reeb orbit $\gamma_{k,\min}$ on $\partial C_k$ and it is nondegenerate. Let $H$ and $H_k$ be Hamiltonians associated to $\partial C$ and $\partial C_k$, respectively. Let $x_{k,\min}$ be the orbits of $H_k$ corresponding to $\gamma_{k,\min}$ via reparametrization. As $k\to\infty$, after passing to a subsequence, $x_{k,\min}$ converges to a periodic orbit of $H$, denoted by $x_{\min}$. Since $\A_{\min}$ is continuous on the space of convex domains in the Hausdorff distance, $\A(x_{\min})= \A_{\min}(\partial C)$.
		
		The proof in Section \ref{sec:pf_A} would give a Floer cylinder with a slit and possibly with boundary punctures $u_k$ in $\R\times\partial C_k$ with $\ev_+(u_k)=x_{k,\min}^-$. A suitable parametrization of $u_k$ converges to a Floer cylinder with a slit and possibly with boundary punctures $u$ in $\R\times\partial C$ with $\ev_+(u)=x_{\min}^-$. Moreover, either (i) or (ii) (see above Lemma \ref{lem:u_nontrivial}) holds for $u$. As before, it suffices to show $E_{d\alpha}(u)>0$. Assume for contradiction that $E_{d\alpha}(u)=0$. Then, as observed in the proof of Lemma \ref{lem:domain_u}, $\mathrm{im\,}u  \subset\R\times x_{\min}$. Since $\Lambda$ does not intersect any minimal periodic orbit, $u$ does not intersect $\R\times \Lambda$. This contradicts that $u$ maps the boundary of its domain to $\R\times\Lambda$.
	\end{rem}
	
	\begin{rem}\label{rem:arbitrary}
	In this remark, we discuss the difficulties in generalizing 	Theorem \ref{thm:chord} to convex domains $C\subset\R^{2n}$, not necessarily strictly convex, that satisfy the Morse-Bott condition. We approximate $C$ by $(C_k)_{k\in\N}$ as in Remark \ref{rem:Lambda}. Then, by the proof in Section \ref{sec:pf_A}, there is a sequence $(u_k)_{k\in\N}$ of Floer cylinders with a slit and possibly with negative boundary punctures in $\R\times \partial C_k$. Now we want to deduce that the limit $u$ of a suitable reparametrization of $u_k$ has positive $d\alpha$-energy as in the proof in this section. 
	
	In the above proof, the asymptotic property of $u$ in \eqref{eq:asymptotic_condition} is essential to obtain $E_{d\alpha}(u)>0$. Note that \eqref{eq:asymptotic_condition} arises from the fact that the component $N$ is of Morse-Bott type and that we perturb $N$ in a specific way, namely by using a Morse function on $\overline{N}$. Since we use this particular perturbation and require $C_k$ to be strictly convex, we have to assume $C$ is also strictly convex, rather than merely convex. 
	
	We do not know whether it is possible to choose an approximation $(C_k)_{k\in\N}$ of an arbitrary convex domain $C$ in such a way that the limit curve $u$ of $u_k$ satisfies an asymptotic property like \eqref{eq:asymptotic_condition}, which would allow us to derive $E_{d\alpha}(u)>0$. It is also unclear to us whether  \eqref{eq:asymptotic_condition} can be obtained for strictly convex domains $C$ without assuming the Morse-Bott condition.
	\end{rem}

\subsection{Proof of Proposition \ref{prop:SH}}

As mentioned in the introduction, Proposition \ref{prop:SH} should follow from a result in \cite{Zho20}. In this section, we outline an alternative proof using the closed-open map in Floer theory, see also the discussion in Remark \ref{rem:improvement} below.

Let $(W,\lambda)$ be an exact symplectic filling of $(Y,\alpha)$, or equivalently $W$ be a Liouville domain with $\partial W=Y$ and $\lambda|_Y=\alpha$, such that its symplectic homology $\SH_*(W)$ vanishes.
Let $\Lambda$ be a closed Legendrian submanifold of $\p W$. Note that as we want to show the non-strict inequality $\A_{\min}(\Lambda)\leq c_\SH(W)$, it suffices to prove the proposition in the case that every periodic Reeb orbit on $Y$ is nondegenerate. Indeed, if this   inequality holds for a contact form approximating $\alpha$, then the inequality holds for $\alpha$ due to the Arzel\'a-Ascoli theorem and the continuity of $c_\SH$. 
By Lemma \ref{lem:c_SH}, we have $c_\SH(W)=c_\RFH(\partial W)$, and it is finite by the hypothesis $\SH_*(W)=0$.  
Assume for contradiction that $c_\SH(W)<\A_{\min}(\Lambda)$. We then take $a\in (c_\RFH(\p W),\A_{\min}(\Lambda))$. We also fix $\epsilon\in (0,\A_{\min}(\p W))$. 

Now we consider the Lagrangian Rabinowitz Floer homology $\RFH^{[0,a)}(\Lambda)$, which is defined using the action functional $\AA_H^\Lambda$ in \eqref{eq:A_lag}, in a manner analogous to the construction of $\RFH^{[0,a)}(\partial W)$, see \cite[Section 8.3]{CO18} and also \cite{Mer14}. This homology is generated by chords with action lower than $a$ and critical points of a Morse function $f_\Lambda$ on $\Lambda$. The boundary operator is defined by counting Floer strips with gradient flow lines of $f_{\Lambda}$. However, since we assumed that there is no Reeb chord with action lower than $a$, the critical points of $f_\Lambda$ are all the generators of $\RFH^{[0,a)}(\Lambda)$, and the boundary operator coincides with the Morse boundary operator for $f_\Lambda$. Hence, $\RFH^{[0,a)}(\Lambda)$ is isomorphic to the singular homology of $\Lambda$. Note that since no Floer strips are involved in the construction, $\RFH^{[0,a)}(\Lambda)$ is well-defined without a Lagrangian filling of $\Lambda$. There is also a map $\iota_{\epsilon,a}^\Lambda: \RFH^{[0,\epsilon)}(\Lambda)\to \RFH^{[0,a)}(\Lambda) $ induced by canonical inclusions at the chain level, cf.~\eqref{eq:iota}. In our case, this is an identity map as there is no generator with action in $[\epsilon,a)$.

We have the following commutative diagram: 
\begin{equation}\label{eq:diagram2}
	\begin{tikzcd}[row sep=1.5em,column sep=1.3em]
	\H_{*-n}(\Lambda;\Z_2) \arrow{r}{\cong} & \RFH^{[0,\epsilon)}(\Lambda)   \arrow{r}{\cong}[swap]{\iota_{\epsilon,a}^\Lambda} &  \RFH^{[0,a)}(\Lambda)    \\
	\H_{*}(\partial W;\Z_2) \arrow{u}{\j} \arrow{r}{\cong} & \RFH^{[0,\epsilon)}(\partial W)    \arrow{u}{\Psi^\epsilon}   \arrow{r}{ \iota_{\epsilon,a} } &  \RFH^{[0,a)}(\partial W) \arrow{u}{\Psi^a}\,.
	\end{tikzcd}	
\end{equation}
The map $\Psi^a$ is the closed-open map defined by counting Floer cylinders with a slit, see the references cited in Section \ref{sec:slit} for a detailed account on this map. We claim that although $\Lambda$ may not have a Lagrangian filling,  $\Psi^a$ is well-defined for a sufficiently stretched almost complex structure due to our assumption. Let $\{J_\nu\}_{\nu\in\N}$ be a sequence of $S^1$-families of almost complex structures on $\widehat W$ stretching the neck along $\{r_*\}\times \partial W$ with $r_* <\log 1/2$, see for example \cite[Section 2.5]{CM18}. Assume for contradiction that there is a sequence of Floer cylinders with a slit $u_\nu:(\Sigma,\partial\Sigma)\to(\R\x\partial C,\R\x\Lambda)$ such that  $\A(\ev_+(u_\nu))<a$ and 
\begin{equation}\label{eq:neck_stretch}
\inf\pi_\R\circ u_\nu|_{\partial\Sigma}<r_*\qquad 
\end{equation}
for sufficiently large $\nu$, where $\pi_\R:\R\x\partial C\to\R$ is the projection along the second factor. Then a suitable reparametrization of $u_\nu$ converges in the $C^\infty_{\rm loc}$-topology to a Floer cylinder with a slit and with negative punctures $u:\Sigma\setminus\Gamma\to\R\times\partial C$ such that $\Gamma$ contains a boundary puncture due to \eqref{eq:neck_stretch}. Thus, $u$ converges at this puncture to a Reeb chord with action lower than $a$. This contradiction proves the claim. 

The map $\j$ is induced by the inclusion $\Lambda\hookrightarrow \p W$ and the Poincar\'e duality. It is isomorphic to closed-open map $\Psi^\epsilon$ for the action window $[0,\epsilon)$.

Let $1_{\p W}$ and $1_\Lambda$ be the fundamental classes in $\H_*(\p W;\Z_2)$ and $\H_*(\Lambda;\Z_2)$, respectively. Since $a>c_\RFH(\p W)$, we have $\iota_{\epsilon,a}(1_{\p W})=0$. This diagram \eqref{eq:neck_stretch} implies $\j(1_{\p W})=1_\Lambda=0$, which is absurd. This contradiction completes the proof.
\qed

\begin{rem}\label{rem:improvement}
The proof of Theorem \ref{thm:chord} carries over to the setting of Proposition \ref{prop:SH}. However, we can only conclude  $\A_{\min}(\Lambda)\leq c_\SH(W)$ with this method as well, as opposed to $\A_{\min}(\Lambda)<c_\SH(W)$, since $\SH_*(W)=0$ only guarantees the existence of a Floer cylinder $v$ with $\ev_+(v)\in S_x$ with $\A(x)=c_{\SH}(W)$, as opposed to $\ev_+(v)=x^-$, see Remark \ref{rem:dyn_convex}. Note that achieving the strict action bound $\A_{\min}(\Lambda)< c_\SH(W)$  anyway does not necessarily yield a Reeb chord with distinct endpoints since $c_\SH(W)$ and $\A_{\min}(\p W)$ do not coincide in general. For this reason, we decided to provide a simpler proof of Proposition \ref{prop:SH} which can be thought of as a homological version of the proof of Theorem \ref{thm:chord}. We finally remark that arguments used in the proof of Proposition \ref{prop:SH} are not enough to deduce the strict action bound  $\A_{\min}(\Lambda)< \A_{\min}(\p C)$ in Theorem \ref{thm:chord}. 
\end{rem}
\begin{rem}
Using the closed-open map in the study of Reeb chords appeared before in \cite{BM14,Rit13,KKK22} in the situation that exact Lagrangian fillings exist. In contrast to those earlier work using the closed-open map for the symplectic homology of $W$, we use  the closed-open map for the Rabinowitz Floer homology of $\p W$.\end{rem}


\bibliographystyle{amsalpha}
\bibliography{chord.bib}
\end{document}